\documentclass
[
    a4paper,
    pdftex,
    onecolumn,
    DIV=11,
    abstracton
]
{scrartcl}

\usepackage
{
    graphicx,
    amssymb,
    amsmath,
    amsthm,
    xcolor,
    dsfont,
    algpseudocode,
    authblk,
    enumerate
}

\usepackage[bf]{caption}
\usepackage[colorlinks,
           pdffitwindow=false,
           plainpages=false,
           pdfpagelabels=true,
           pdfpagemode=UseOutlines,
           pdfpagelayout=SinglePage,
           bookmarks=false,
           colorlinks=true,
           hyperfootnotes=false,
           linkcolor=blue,
           citecolor=green!50!black]
{hyperref}

\DeclareMathAlphabet{\mathpzc}{OT1}{pzc}{m}{it}

\newcommand{\subfiguretitle}[1]{{\scriptsize{#1}} \\[1mm] }
\newcommand{\R}{\mathbb{R}}
\newcommand{\C}{\mathbb{C}}
\newcommand{\N}{\mathbb{N}}
\newcommand{\cB}{{\mathcal{B}}}
\newcommand{\cC}{{\mathcal{C}}}
\newcommand{\cM}{{\mathcal{M}}}
\newcommand{\cN}{{\mathcal{N}}}
\newcommand{\cU}{{\mathcal{U}}}

\DeclareMathOperator{\diam}{diam}

\newtheorem{theorem}{Theorem}[section]

\newtheorem{proposition}[theorem]{Proposition}
\newtheorem{definition}[theorem]{Definition}
\newtheorem{remark}[theorem]{Remark}
\newtheorem{example}[theorem]{Example}
\newtheorem{algorithm}[theorem]{Algorithm}

\makeatletter
\renewcommand*\env@matrix[1][*\c@MaxMatrixCols c]{%
  \hskip -\arraycolsep
  \let\@ifnextchar\new@ifnextchar
  \array{#1}}
\makeatother

\begin{document}

\title{A Set-Oriented Numerical Approach for Dynamical Systems with Parameter Uncertainty}
\author[1]{Michael Dellnitz}
\author[2]{Stefan Klus}
\author[1]{Adrian Ziessler}
\affil[1]{\normalsize Department of Mathematics, Paderborn University, Germany}
\affil[2]{\normalsize Department of Mathematics and Computer Science, Freie Universit\"at Berlin, Germany}


\maketitle

\begin{abstract}
In this article, we develop a set-oriented numerical methodology which allows
to perform uncertainty quantification (UQ) for dynamical systems from a
global point of view. That is, for systems with uncertain parameters we
approximate the corresponding global attractors and invariant measures
in the related stochastic setting. Our methods do not rely on generalized 
polynomial chaos techniques. Rather, we extend classical set-oriented
methods designed for deterministic dynamical systems \cite{DH97, DJ99}
to the UQ-context,
and this allows us to analyze the long-term uncertainty propagation. The algorithms
have been integrated into the software package GAIO \cite{DFJ01}, and
we illustrate the use and efficiency of these techniques by a couple of numerical examples.
\end{abstract}


\section{Introduction}

The analysis of the influence of uncertainties in complex dynamical systems on the system's behavior has gained considerable attention in the last years. In many applications, input parameters, initial conditions, or boundary conditions are not known exactly and are thus described by probability distributions. The goal is to quantify the effects of these uncertainties and their impact on, for instance, stability or performance of the system. Due to the stretching and folding of the corresponding trajectories, propagating probability density functions through highly nonlinear dynamical systems is particularly challenging~\cite{LBR14}.

In this article, we consider parameter-dependent discrete dynamical systems assuming that some parameters are uncertain. The main goal is to develop robust algorithms for the analysis of the resulting statistical behavior of the
system so that we capture the long-term uncertainty propagation. To this end, we compute approximations of the corresponding \emph{invariant sets} and \emph{invariant measures} using so-called \emph{set-oriented numerical methods}. These have been developed for the numerical analysis of complex dynamical systems, see e.g.~\cite{DH97, DJ99, FD03, FLS10}, and used for a host of different application areas such as molecular dynamics~\cite{SHD01}, astrodynamics~\cite{DJLMPPRT05}, and ocean dynamics~\cite{FHRSSG12}. Recently, set-oriented methods have been extended to compute attractors for delay differential differential equations~\cite{DHZ16}. The basic idea is to cover the objects of interest by outer approximations created via multilevel subdivision schemes. In this
work, we will generalize these techniques to the context of uncertainty quantification. 

A related -- though not set-oriented -- approach using generalized polynomial chaos (gPC),
see e.g.~\cite{Sud08, Sul15}, has been utilized in~\cite{LBR14} for
the numerical analysis of uncertainty in dynamical systems. Long-term uncertainty propagation is accomplished by approximating and composing intermediate short-term flow maps using spectral polynomial bases. A gPC based approach has also been applied in \cite{PS14} in the context of Hamiltonian systems exhibiting multi-scale
dynamics, or in \cite{Xiu07} with an application to differential algebraic equations (DAEs).
In contrast to this, our work relies on efficient adaptive subdivision schemes and transfer-operator based methods. This is the first time a set-oriented approach has been used for the computation of attractors and invariant measures in this context,
and this allows us to quantify the uncertainty from a global point of view. Technically, the results of this article are mainly based on two previous publications: First, we extend the classical subdivision scheme in \cite{DH97} to the context of parameter dependence within a compact set $\Lambda$. Secondly, our measure computations in Section~\ref{sec:comp_para_meas} are based on the theoretical framework in~\cite{DJ99}.

A detailed outline of the paper is as follows: In Section~\ref{sec:comp_para_att}, we introduce the notion of \emph{$(Q,\Lambda)$-attractors}. This is an object in state space which contains all the invariant sets within the compact set $Q$ which can potentially be created by the parameter uncertainty within $\Lambda$. That is, the particular probability distribution on $\Lambda$ is not yet relevant for our considerations in this section. Furthermore, we develop an algorithm which allows to compute outer approximations of $(Q,\Lambda)$-attractors (Algorithm~\ref{alg:subdivision}). In Section~\ref{sec:comp_para_meas}, we assume that the uncertainty within $\Lambda$ is given by a certain probability distribution. We model this parameter uncertainty with appropriate {\em stochastic transition functions} and develop an algorithm for the computation of corresponding invariant measures (Algorithm~\ref{alg2}). The use of {\em small random perturbations} \cite{K86} allows us to prove a related convergence result, which is essentially an adapted extension of the corresponding result in \cite{DJ99}. In Section \ref{sec:num_exam}, we illustrate the use and efficiency of the algorithms by three examples, namely the H\'enon map, the van der Pol oscillator and the Arneodo system. Finally, in Section~\ref{sec:Conclusion}, we conclude with a short summary of the main results and possible future work.

\section{Computation of $(Q,\Lambda)$-Attractors via Subdivision}
\label{sec:comp_para_att}

In this section, we will introduce the notion of $(Q,\Lambda)$-attractors, a generalization of the concept of relative global attractors defined in \cite{DH97}. The goal is to develop a numerical technique which allows
to identify the region in state space which is potentially influenced by inherent parameter uncertainties.

\subsection{$(Q,\Lambda)$-Attractors}
\label{ssec:QLA}

Let us denote by $\Lambda \subset \R^p$ a compact subset which represents our set of admissible parameter values. Then for each $\lambda \in \Lambda$ we have the dynamical system
\begin{equation} \label{eq:dynsys1}
    x_{j+1} = f(x_j,\lambda), \quad j = 0, 1, \dots,
\end{equation} 
where $x_j\in \R^n$ and $f: \R^n \times \Lambda \to \R^n$ is continuous. In order to take the uncertainty in \eqref{eq:dynsys1} with respect to $\lambda$ into account, we consider the corresponding (set-valued) map $F_\Lambda : \R^n \to {\mathcal P}(\R^n)$, where
\begin{equation*}
    F_\Lambda(x) = f(x,\Lambda)
\end{equation*} 
and $ {\mathcal P}(\R^n)$ denotes the power set of $\R^n$.
 
The purpose of this section is to approximate the region in state space covering all the backward invariant sets which can potentially be generated by $\lambda$-distributions with support in $\Lambda$. Therefore, it is not yet necessary to assume that the uncertainty with respect to the parameters is described by a certain distribution. However, we will come back to this in Section~\ref{sec:comp_para_meas} where we compute invariant measures for different $\lambda$-distributions.

In this section, we develop an algorithm which allows to compute the object
\begin{equation} \label{eq:relativeAttractor}
    A_{Q,\Lambda} = \bigcap_{j\ge 0} F_\Lambda^j(Q),
\end{equation}
where $Q \subset \R^n$ is a compact subset. We call $A_{Q,\Lambda}$ the \emph{$(Q,\Lambda)$-attractor}.

\begin{remark} \label{rem:classic}
\begin{enumerate}[(a)]
\item With \eqref{eq:relativeAttractor}, we generalize the concept of \emph{relative global attractors} which has first been introduced in \cite{DH97} for dynamical systems of the form
\begin{equation*}
    x_{j+1} = g(x_j), \quad j= 0,1,\dots,
\end{equation*}
where $g$ is a homeomorphism. In fact, the \emph{global attractor relative to} $Q$ is defined by
\begin{equation*}
    A_Q = \bigcap_{j\ge 0} g^j(Q).
\end{equation*}
\item In \cite{DSchS02,DSchH05}, relative global attractors have been defined for non-autonomous dynamical systems. More precisely, given $s$ dynamical systems $g_1, \dots, g_s$, one is interested in computing their common invariant sets. Therefore, the object of interest is
\begin{equation*}
    A_{Q,g_1,\dots,g_s} = \bigcap_{\omega\in\Omega} \bigcap_{j \ge 1} g_{\omega^j}(Q)\cap Q.
\end{equation*}
Here, $\Omega = \{1,2,\dots,s\}^{\N_0}$ is the space of sequences of $s$ symbols, and for $\omega=(\omega_i)\in \Omega$ the map $g_{\omega^j}$ is defined as the composition $g_{\omega^j} = g_{\omega_{j-1}}\circ \cdots\circ g_{\omega_0}$.

\item Our definition in \eqref{eq:relativeAttractor} is also strongly related to concepts which have been introduced in connection with control problems. For instance, the (positive and negative) \emph{viability kernels} as defined in \cite{S01} are precisely of this type. However, our analytical setup is different and it is the purpose of this article to utilize the related set-oriented numerical analysis in the context of uncertainty quantification.
\end{enumerate}
\end{remark}

The following properties of $A_{Q,\Lambda}$ follow immediately from its definition.

\begin{proposition} \label{prop:AQprops}
\begin{enumerate}[(a)]
\item Suppose that the map $F_\Lambda : \R^n \to {\mathcal P}(\R^n)$ is one-to-one. Then $A_{Q,\Lambda}$ is backward invariant, i.e.\
\begin{equation*}
    F_\Lambda^{-1}(A_{Q,\Lambda}) \subset A_{Q,\Lambda}.
\end{equation*}
\item Let $B \subset Q$ be a set such that $B\subset F_\Lambda(B)$. Then $B\subset A_{Q,\Lambda}$.
\end{enumerate}
\end{proposition}

By this proposition, $A_{Q,\Lambda}$ should be viewed as the set which contains all the backward invariant sets of $F_\Lambda$. Accordingly, in the deterministic context $A_Q$ typically consists of all invariant sets and their unstable manifolds within $Q$.

\subsection{Subdivision Scheme}
\label{ssec:subdiv}

The numerical developments in this work are based on a modification of the subdivision scheme introduced in \cite{DH97}, where this algorithm is used for an approximation of relative global attractors (see Remark~\ref{rem:classic} (a)). Here, we extend it in order to create coverings of $(Q,\Lambda)$-attractors.

This algorithm generates a sequence $\cC_0,\cC_1,\dots$ of finite collections of compact subsets of $\R^n$ such that the diameter
\begin{equation*}
    \diam(\cC_\ell) = \max_{C\in\cC_\ell}\diam(C)
\end{equation*}
converges to zero for $\ell\rightarrow\infty$. Concretely, given an initial collection $\cC_0$ with $\bigcup_{C\in\cC_0}C = Q$, we successively obtain $\cC_\ell$ from $\cC_{\ell-1}$ for $\ell=1,2,\dots$ in two steps:

\begin{enumerate}
\item \emph{Subdivision:} Construct a new collection $\hat\cC_\ell$ such that
\begin{equation} \label{eq:sd1C}
    \bigcup_{C\in\hat\cC_\ell}C = \bigcup_{C\in\cC_{\ell-1}}C
\end{equation}
and
\begin{equation} \label{eq:sd2C}
    \diam(\hat\cC_\ell) = \theta_\ell\diam(\cC_{\ell-1}),
\end{equation}
where $0<\theta_{\min} \le \theta_\ell\le \theta_{\max} < 1$.
\item \emph{Selection:} Define the new collection $\cC_\ell$ by
\begin{equation} \label{eq:selectC}
    \cC_\ell=\left\{C\in\hat\cC_\ell : \exists \hat C\in\hat\cC_\ell \text{ and } \exists
     \lambda\in \Lambda \text{ such that } f(\cdot,\lambda)^{-1}(C)\cap\hat C \ne \varnothing \right\}.
  \end{equation}
\end{enumerate}

We denote by $Q_\ell$ the area in state space covered by $\cC_\ell$, that is
\begin{equation*} 
    Q_\ell=\bigcup_{C\in\cC_\ell}C,
\end{equation*}
and let
\begin{equation*}
    Q_\infty = \lim\limits_{\ell \rightarrow\infty} Q_\ell.
\end{equation*}
Observe that this limit exists due to the fact that the $Q_\ell$ form a sequence of nested compact sets.

\begin{remark} \label{rem:classic1}
\begin{enumerate}[(a)]
\item If we replace $f(x,\lambda)$ by $g=g(x)$ in the subdivision scheme ($g$ a homeomorphism) 
then the main theoretical result in \cite{DH97} states that 
\begin{equation*}
    A_Q = Q_\infty,
\end{equation*}
see also Remark~\ref{rem:classic} (a).
\item The result in (a) is also valid in the situation where $g$ is just continuous -- and not a homeomorphism. However, in this case one has to assume additionally that $g^{-1}(A_Q) \subset A_Q$. This result has recently been proved in \cite{DHZ16}.
\end{enumerate}
\end{remark}

Now we state the main result of this section:
\begin{proposition}
Suppose that the map $F_\Lambda : \R^n \to {\mathcal P}(\R^n)$ is one-to-one. Then
\begin{equation*}
    A_{Q,\Lambda} = Q_\infty.
\end{equation*}
\end{proposition}

The proof is in principle identical to the one in \cite{DH97} or \cite{DHZ16}, respectively. Thus, we just sketch it here.

\begin{proof}[Sketch of Proof]
The proof consists of the following steps:
\begin{enumerate}[(i)]
\item $A_{Q,\Lambda} \subset Q_\ell$ for all $\ell$. Here, one uses Proposition~\ref{prop:AQprops} (a).
\item $Q_\infty \subset F_\Lambda (Q_\infty)$. In this step, the continuity of $f$ is crucial.
\end{enumerate}
Now by (i)
\begin{equation*}
    A_{Q,\Lambda} \subset Q_\infty,
\end{equation*}
and by (ii) and Proposition~\ref{prop:AQprops} (b)
\begin{equation*}
    Q_\infty \subset A_{Q,\Lambda},
\end{equation*}
yielding the desired result.
\end{proof}

Generalizing the approach in \cite{DH97}, we use the following numerical realization of the
subdivision scheme for the approximation of the $(Q,\Lambda)$-attractor.

\begin{algorithm}\label{alg:subdivision}
    Choose an initial box $Q \subset \R^n$, defined by a generalized rectangle of
    the form
    \begin{equation*}
        Q(c,r) = \{ y\in \R^n : |y_i-c_i|\leq r_i \text{ for } i = 1,\dots,n \},
    \end{equation*}
    where $c,r \in \R^n,\ r_i > 0$ for $i=1,\ldots,n$, are the center and the radius, 
    respectively. Discretize $\Lambda = \{\lambda_1, \dots, \lambda_M\}$ uniformly and start the subdivision
    algorithm with a single box $\cC_0 = \{ Q \}$.
    \begin{enumerate}
        \item Realization of the subdivision step: In step $(\ell-1)$, we
        subdivide each box $C \in \cC_{\ell-1}$ of the current collection by bisection with respect
        to the $s$-th coordinate, where $s$ is varied cyclically. Thus, in the new collection
        $\hat{\cC}_{\ell}$ the number of boxes is increased by a factor of $2$ (cf.~\eqref{eq:sd1C}, \eqref{eq:sd2C}).
        \item Realization of the selection step: We choose a finite set of test points in each box
        $C_j \in \cC_\ell$ and replace the condition \eqref{eq:selectC} by
        \begin{align}\label{eq:selectCmod}
            f(x,\lambda_k) \notin C_i \text{ for all test points } x\in C_j \text{ and all } \lambda_k,\ k = 1,\dots,M.
        \end{align}A box $C_i$ is discarded if \eqref{eq:selectCmod} is satisfied for all $j$.
        \item Repeat (1)+(2) until a prescribed size $\varepsilon$ of the diameter relative to the initial box $Q$ is reached. That is, we stop when
        \begin{equation*}
            \diam(\cC_\ell) < \varepsilon \diam(Q).
        \end{equation*}
    \end{enumerate}
\end{algorithm}

\begin{remark}
\begin{enumerate}[(a)]
\item The distribution of test points in each box  $C \in \cC_\ell$ is performed as follows.
Observe that each box is defined by a generalized rectangle with a radius $r$ and center $c$ and therefore is the affine image of the standard cube $[-1,1]^n$ scaled by $r$ and translated by $c$. Thus, by using this transformation it is sufficient to define the distribution of test points for the standard cube, e.g.~via a uniform grid or a (quasi-)Monte Carlo sampling. In our computation, we use the Halton sequence, which is a quasi-random number sequence~\cite{H64algorithm}.
\item Algorithm~\ref{alg:subdivision} has been developed within the software package GAIO
(see \cite{DFJ01}). Thus, as in the deterministic case the numerical complexity depends crucially
on the dimension of $A_{Q,\Lambda}$. In fact, our experience indicates that we can approximate 
these objects for dimensions up to three even in higher dimensional state space. On the other hand,
it will become very time-consuming to compute  $A_{Q,\Lambda}$ if its dimension is larger than four.
\end{enumerate}
\end{remark}
\section{Computation of Invariant Measures on $(Q,\Lambda)$-Attractors}
\label{sec:comp_para_meas}

We use a transfer operator approach to approximate invariant measures on $(Q,\Lambda)$-attractors.
This is a classical mathematical tool for the numerical analysis of complicated dynamical
behavior, e.g.\ \cite{DJ99, SHD01, FP09, koltai2011}, and here we use it in the context of
uncertainty quantification.

\subsection{Invariant Measures and the Transfer Operator}

We briefly introduce the reader to the notion of \emph{transfer operator} in the stochastic setting.
In the context of stochastic differential equations, transfer operators and stochastic
transition functions have recently been utilized in \cite{FK15}. Also in this paper, we have to stay within the stochastic setting in order to take the parameter uncertainty into account. In the following paragraphs, we follow closely the related contents in \cite{DJ99}.

As in the previous section, let $Q\subset \R^n$ be compact and denote by $\cB$ the Borel-$\sigma$ algebra on $Q$.
\begin{definition}\label{def:stf}
A function $p:Q\times\cB\to\R$ is a \emph{stochastic transition function}, if 
\begin{enumerate}
\item $p(x,\cdot)$ is a probability measure for every $x\in Q$,
\item $p(\cdot, A)$ is Lebesgue-measurable for every $A\in\cB$.
\end{enumerate}
\end{definition}

\begin{remark}\label{rem:deter}
Let $\delta_y$ denote the Dirac measure supported on the point $y\in Q$. Then $p_{\bar\lambda}(x,A)=\delta_{f(x,\bar\lambda)}(A)$, $\bar\lambda \in \Lambda$, is a stochastic transition function. This represents the deterministic situation for fixed $\bar\lambda\in\Lambda$ in this stochastic setup.
\end{remark}

We now define the notion of an invariant measure in the stochastic setting. To this end, we denote by $\cM$ the set of probability measures on $\cB$.
\begin{definition}
Let $p$ be a stochastic transition function. If $\mu\in\cM$ satisfies
\begin{equation*}
    \mu(A) = \int p(x,A)\;d\mu(x)
\end{equation*}
for all $A\in\cB$, then $\mu$ is an \emph{invariant measure of $p$}.
\end{definition}

The following example illustrates the previous remark that we recover the deterministic situation in the case where $p_{\bar\lambda}(x,\cdot)=\delta_{f(x,\bar\lambda)}$.
\begin{example} \label{ex:invclassic}
Suppose that $p_{\bar\lambda}(x,\cdot)=\delta_{f(x,\bar\lambda)}$ and let $\mu$ be an invariant measure of $p_{\bar\lambda}$.  Then we compute for $A\in\cB$
\begin{eqnarray*}
    \mu(A) &=& \int p_{\bar\lambda}(x,A)\;d\mu(x) =
    \int \delta_{f(x,\bar\lambda)}(A) \;d\mu(x)\\ &=& \int \chi_A(f(x,\bar\lambda)) \;d\mu(x)
    = \mu(f^{-1}(A,\bar\lambda)),
\end{eqnarray*}
where we denote by $\chi_A$ the characteristic function of $A$. Hence, $\mu$ is an invariant measure for the map $f(\cdot,\bar\lambda)$ in the classical deterministic sense (cf.\ \cite{Po93}).
\end{example}


\begin{definition}
Let $p$ be a stochastic transition function.  Then the \emph{transfer operator\/} $P:\cM_\C\to\cM_\C$ is defined by
\begin{equation*}
    P\mu(A) = \int p(x,A)\;d\mu(x),
\end{equation*}
where $\cM_\C$ is the space of bounded complex-valued measures on ${\cal B}$.
\end{definition}

By definition, an invariant measure $\mu$ is a fixed point of $P$, i.e.
\begin{equation} \label{eq:pmu}
    P \mu = \mu,
\end{equation}
and in the remainder of this section we develop a numerical method for the approximation of such measures.

\subsection{Stochastic Transition Functions on $(Q,\Lambda)$-Attractors}

Suppose that the parameter uncertainty on $\Lambda$ is given by the probability density function $\rho:\Lambda \to \R_{\ge 0}$. Then we define the corresponding stochastic transition function $q_\rho$ as follows: For $A\subset \R^n$ and each point $x\in \R^n$ let
\begin{equation*}
    \Lambda_x(A) = \left\{ \lambda \in \Lambda : f(x,\lambda) \in A \right\}= f(x,\cdot)^{-1}(A).
\end{equation*}
The set $\Lambda_x(A)$ is measurable for each $x\in \R^n$ and each $A\in \cB$ since $f$ is continuous, and we define the measure of $A$ to be the measure of $\Lambda_x(A) \subset \Lambda$, that is
\begin{equation*}
    \label{eq:p}
    q_\rho (x,A) = \int\limits_{\Lambda_x(A)} \rho(\lambda) \; d\lambda.
\end{equation*}
Observe that
\begin{equation}\label{eq:density}
    q_\rho (x,A) = \int\limits_{f(x,\cdot)^{-1}(A)} \rho(\lambda) \; d\lambda =
    \int\limits_{\Lambda} \chi_A(f(x,\lambda)) \rho(\lambda) \; d\lambda.
\end{equation}
By construction, $q_\rho (x,\cdot)$ is a probability measure for every $x\in\R^n$.
Moreover, $q_\rho(\cdot,A)$ is integrable by (\ref{eq:density}) and therefore a stochastic transition function according to Definition~\ref{def:stf}.

In the particular deterministic case where the parameter $\bar\lambda$ is fixed, we obtain (see Remark~\ref{rem:deter})
\begin{equation*}
    q_\rho (x,A) = \int\limits_{\Lambda} \chi_A(f(x,\lambda)) \delta(\bar\lambda) \; d\lambda
    = \chi_A(f(x,\bar\lambda)) = \delta_{f(x,\bar\lambda)}(A) = p_{\bar\lambda}(x,A).
\end{equation*}
Here $\delta(\cdot)$ is the Dirac delta function.

\subsection{Approximation of Invariant Measures}

Given a stochastic transition function $q_\rho (x,\cdot)$, we now explain how to approximate a corresponding invariant measure numerically. For the discretization of the problem, we proceed as in \cite{dellnitz1997exploring,koltai2011} and use characteristic functions $\chi_{C_l}$ $(l=1,2,\ldots,d)$
on each box in our covering of the $(Q,\Lambda)$-attractor which has been generated by Algorithm~\ref{alg:subdivision}. We denote by $m$ the Lebesgue measure and define the corresponding probability measures
\begin{equation*}
    \mu_{C_l} (A) =  \frac{1}{m(C_l)} \displaystyle\int_{A} \chi_{C_l} \; dm, \quad l=1,\dots,d.
\end{equation*}
The transfer operator $P$ is acting on these measures as follows
\begin{equation*}
    \left( P \mu_{C_l}\right) (A) = \int q_\rho(x,A) \; d\mu_{C_l}
    = \frac{1}{m(C_l)}\displaystyle\int q_\rho(x,A) \chi_{C_l}\, dm
    = \frac{1}{m(C_l)}\displaystyle\int_{C_l} q_\rho(x,A) \; dm.
\end{equation*}

Thus, we can approximate the transfer operator with the following stochastic matrix $P_d = \left(p_{kl} \right)$ on the box covering $\{C_1,\dots,C_d\}$:
\begin{equation}\label{eq:pkl}
    p_{kl} = \frac{1}{m(C_l)} \displaystyle\int_{C_l} q_\rho(x,C_k) \; dm, \quad k,l=1,\dots,d.
\end{equation}

Finally, we approximate the invariant measure corresponding to the stochastic transition function $q_\rho$ by the stationary distribution of the Markov chain given by $P_d$. Concretely, we approximate probability measures $\nu\in\cM$ by
\begin{equation*}
    \nu \approx \sum_{l=1}^d \alpha_l \,\mu_{C_l}.
\end{equation*}
In order to obtain an approximation of an invariant measure $\mu$, we require
\begin{equation*}
    \left( P\sum_{l=1}^d \alpha_l \,\mu_{C_l}\right) (C_k) = \sum_{l=1}^d \alpha_l \,\mu_{C_l}(C_k) = \alpha_k,
    \quad k=1,2,\ldots,d.
\end{equation*}
Here, we have used the fact that
\begin{equation*}
    \mu_{C_l} (C_k) = \delta_{kl}\quad \mbox{($\delta_{kl}$ the Kronecker-delta)}
\end{equation*}
by the construction of the box collection. That is, for an approximation of
an invariant measure $\mu$ we have to compute the eigenvector $\alpha_d\in\R_{\ge 0}^d$ for the eigenvalue $ 1 $ of the matrix $P_d$, i.e.
\begin{equation*}
    P_d \,\alpha_d = \alpha_d
\end{equation*}
(see also \eqref{eq:pkl}).

\begin{remark}
\label{rem:to}
\begin{enumerate}[(a)]
\item In the deterministic case, that is $q_\rho (x,\cdot) =
p_\lambda(x,\cdot)=\delta_{f(x,\lambda)}$, the transition probabilities are given by
\begin{equation*} \label{eq:pkl_det}
    p_{kl} = \frac{m\left(f^{-1}(C_k,\lambda)\cap C_l\right)}{m(C_l)}.
\end{equation*}
\item Numerically, the computation of \eqref{eq:pkl} is realized as follows: For each $1 \leq l \leq d$, select test points $x_1,\dots,x_N \in C_l$ and
uncertain parameters $ \lambda_1,\dots,\lambda_M $ 
distributed according to the probability density function $\rho$. This yields
\begin{align*}
    p_{kl} &= \frac{1}{m(C_l)} \displaystyle\int_{C_l} q_\rho(x,C_k) \; dm\\
    &\approx \frac{1}{M\cdot N}\sum_{i=1}^{M} \big\vert \left\{ j \in \{ 1,\ldots,N\} | f(x_j,\lambda_i) \in C_k \right\}\big\vert.
\end{align*}
\item The box covering $\{ C_1,\dots,C_d\}$ obtained with the subdivision scheme and the dynamics induced by the stochastic transition function $q_\rho$ yield a directed graph as illustrated in Figure~\ref{fig:graph}. The dynamics on this graph with the transition probabilities in (\ref{eq:pkl}) can be viewed as an approximation of the transfer operator $P$.

\begin{figure}[htb]
    \centering
    \includegraphics[width=0.8\textwidth]{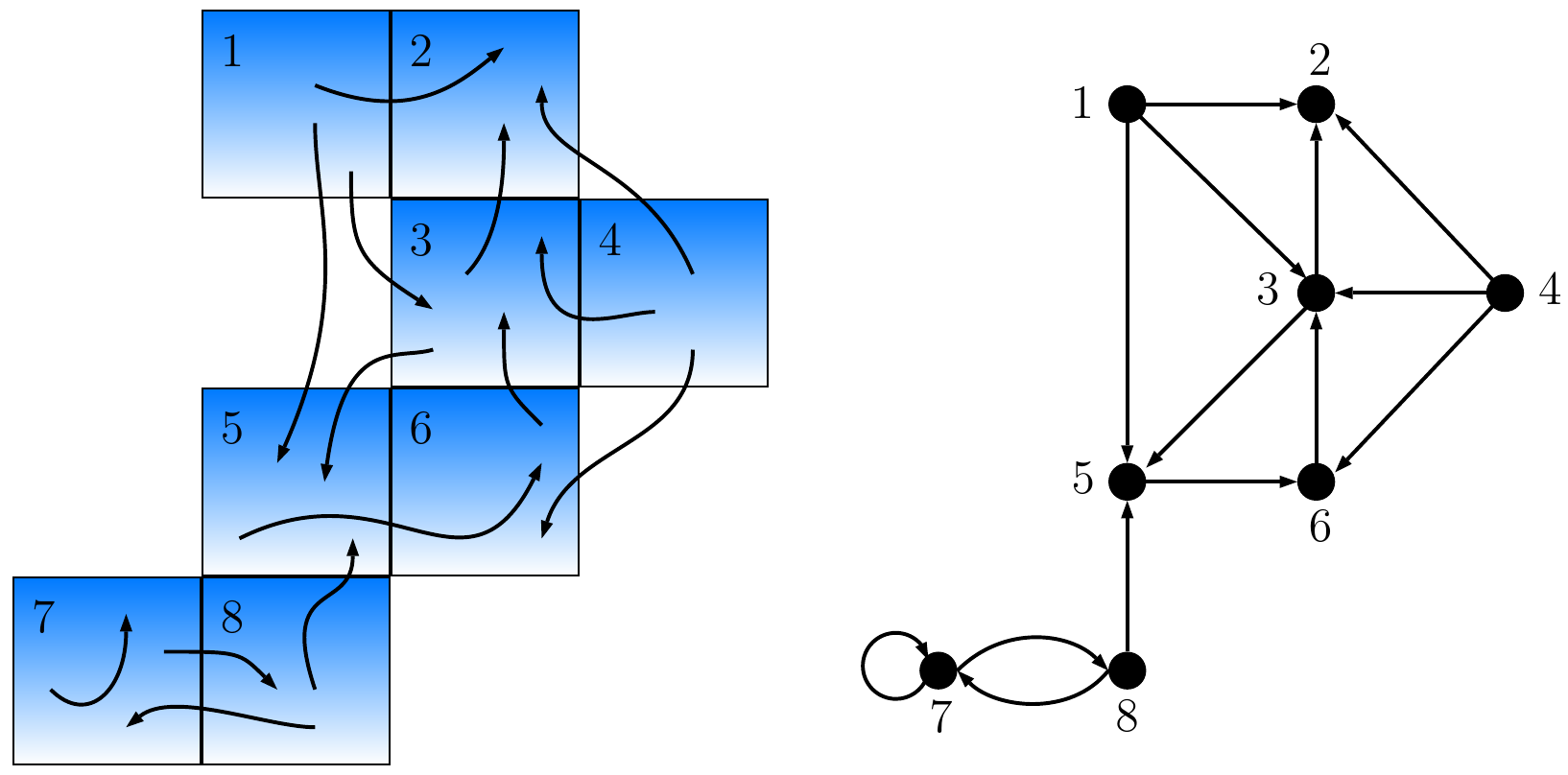}
    \caption{Schematic illustration of the graph induced by the transition function $ q_\rho $ on the box covering. Left: Box covering $ \{C_1, \dots, C_d\} $ and mapping of points from box $ C_l $ to box $ C_k $. Right: Resulting directed graph with vertices $ \{v_1, \dots, v_d\} $ and edges $ (v_l, v_k) $.}
    \label{fig:graph}
\end{figure}

\end{enumerate}
\end{remark}

We summarize our numerical approach in the following algorithm.
\begin{algorithm} \label{alg2}
    The strategy for the approximation of an invariant measure corresponding to the stochastic transition function $q_\rho$ supported on a $(Q,\Lambda)$-attractor $A_{Q,\Lambda}$ can now be formulated as follows:
    \begin{enumerate}
        \item Approximate a $(Q,\Lambda)$-attractor $A_{Q,\Lambda}$ by the subdivision Algorithm~\ref{alg:subdivision} to obtain a box covering $\{C_1,\dots,C_d\}$.
        \item Use $\{C_1,\dots,C_d\}$ to compute the discretized transfer operator $P_d$ by \eqref{eq:pkl}. If $\lambda$ is distributed according to the probability density function $\rho$, we use Remark~\ref{rem:to} (b) where the $\lambda$-values are obtained via a (quasi-)Monte Carlo sampling.        
        \item Compute the eigenvector $\alpha_d \in \R^d$ corresponding to the eigenvalue $1$ of $P_d$ to obtain an approximation of an invariant measure $\mu$ on $A_{Q,\Lambda}$ (cf.~\eqref{eq:pmu}).
    \end{enumerate}
\end{algorithm}

\subsection{Convergence Result}

We utilize the theoretical framework from \cite{DJ99} to obtain a convergence result for our numerical approach. Our developments in this section cover the classical deterministic situation (see Example~\ref{ex:invclassic}). Therefore, we have to consider \emph{small random perturbations} (cf.\ \cite{K86}) of $f(x,\lambda)$ so that the transfer operator becomes compact as an operator on $L^2$. That is, in addition to the inherent parameter uncertainty we now introduce a perturbation in state space so that classical convergence theory for compact operators can be applied.

We let $B=B_0(1)$ the open ball in $\R^n$ of radius one and define for
$\epsilon > 0$ 
\begin{equation*}
    k_\epsilon(x,y)=
    \frac{1}{\epsilon^n m(B)}\chi_B\left(\frac{1}{\epsilon}\Big(y-x\Big)\right),
    \quad x,y\in \R^n.
\end{equation*}
We use this function for the definition of a {\em transition density function} which allows to take the parameter uncertainty into account
\begin{equation*} \label{eq:k-eps}
    k_{\epsilon,f} (x,y) = \int\limits_{\Lambda} k_\epsilon(f(x,\lambda),y) \rho(\lambda) \; d\lambda.
\end{equation*}
With this we define a stochastic transition function $p_\epsilon$ in this
random context by
\begin{equation*} \label{eq:p-eps}
    p_\epsilon(x,A) = \int_A k_{\epsilon,f}(x,y)\;dm(y).
\end{equation*}
\begin{remark}
Observe that
\begin{equation*}
    \int\limits_A k_\epsilon(f(x,\lambda),y) \;dm(y) \to \delta_{f(x,\lambda)}(A) \quad \text{for } \epsilon \to 0,
\end{equation*}
and therefore
\begin{equation*}
    p_\epsilon(x,A) = \int_A k_{\epsilon,f}(x,y)\;dm(y) \to 
    \int\limits_{\Lambda} \chi_A(f(x,\lambda)) \rho(\lambda) \; d\lambda \quad \text{for } \epsilon \to 0.
\end{equation*}
Thus, as expected, we obtain in the limit the stochastic transition function
$q_\rho (x,A)$ in \eqref{eq:density}.
\end{remark}

Due to the small random perturbation, the measure $p_\epsilon(x,\cdot)$ is absolutely continuous for $\epsilon > 0$, and the corresponding transfer operator $P_\epsilon : L^1 \to L^1$ can be written as
\begin{equation} \label{eq:P-eps}
    \left( P_\epsilon g\right) (y) = \int k_{\epsilon,f}(x,y) g(x)\;dm(x)\quad \text{for all } g\in L^1.
\end{equation}

It is easy to verify that for $\epsilon > 0$
\begin{equation*}
    \iint |k_{\epsilon,f}(x,y)|^2\;dm(x)dm(y) < \infty.
\end{equation*}
Therefore, the transfer operator in \eqref{eq:P-eps} as an operator $P_\epsilon:L^2\to L^2$ is compact.

Now let $\beta \not=0$ be an eigenvalue of $P_\epsilon$ and let $E$ be a projection onto the corresponding generalized eigenspace. Then we have the following convergence result which yields an approximation result for invariant measures in the randomized situation (see Theorem 3.5 in \cite{DJ99} and also \cite{osborn75}).

\begin{theorem}
Let $\beta_d$ be an eigenvalue of $P_d$ such that $\beta_d\to\beta$ for $d\to\infty$, and let $\gamma_d$ be a corresponding eigenvector of unit length. Then there is a vector $h_d\in R(E)$ and a constant $C>0$ such that $(\beta I-P)h_d = 0$ and
\begin{equation*}
    \|h_d-\gamma_d\|_2 \le C\|(P_\epsilon-P_d)|_{R(E)}\|_2.
\end{equation*}
\end{theorem}

\section{Numerical Results}
\label{sec:num_exam}
In this section, we will present numerical results for different dynamical systems. For each system, we assume that there exists one uncertain parameter, denoted by $\lambda$. In what follows, let $ \cU(\Lambda) $ denote a uniform distribution over the set $ \Lambda $ and $ \cN(\mu,\sigma^2)$ a Gaussian -- if necessary truncated so that it fits into $\Lambda$ -- with
mean $ \mu $ and standard deviation $ \sigma $.

\subsection{H\'enon Map}
\begin{figure}[htb]
    \centering
    \includegraphics[width=.8\textwidth]{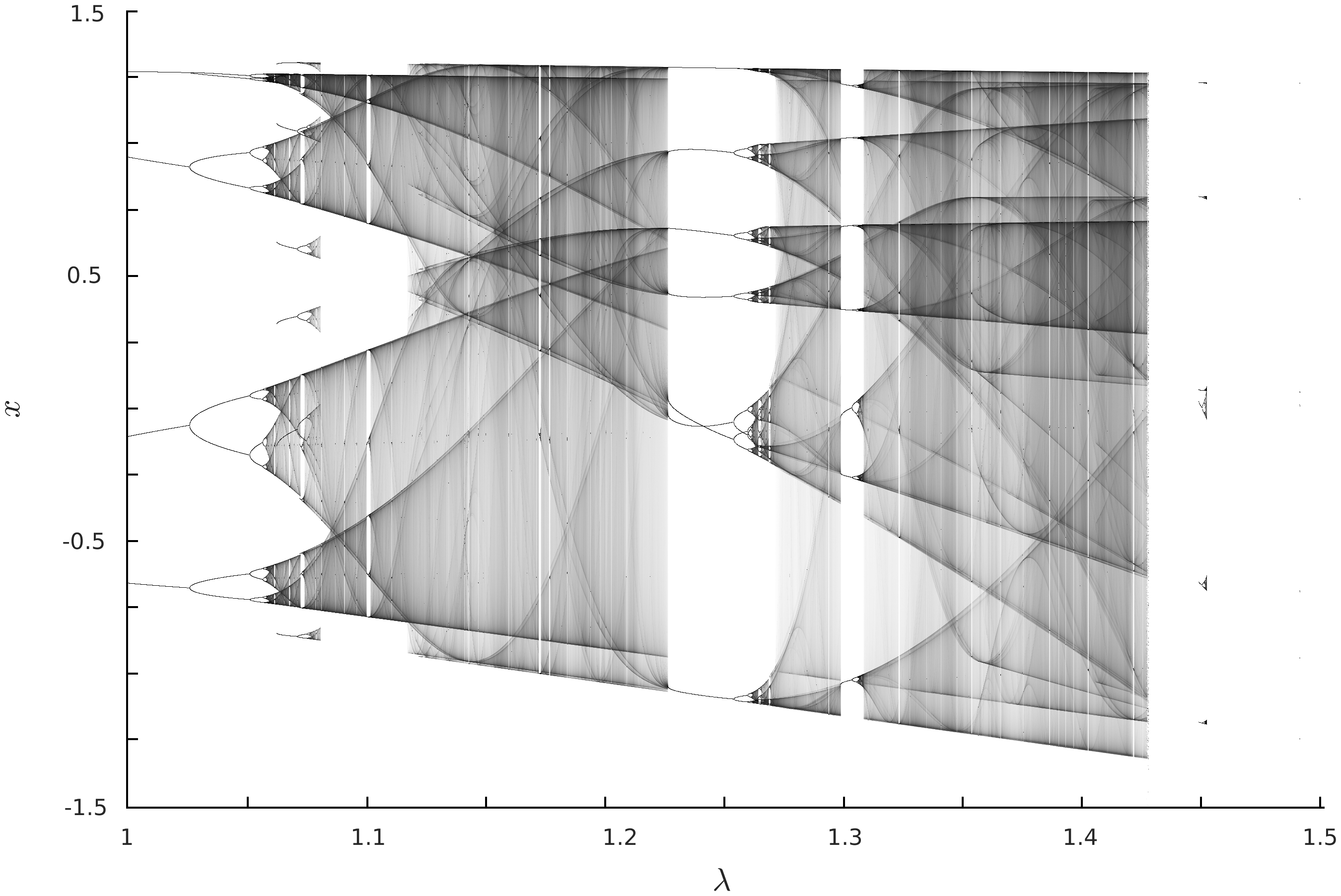}
    \caption{Bifurcation diagram for the H\'enon map with $\nu = 0.3$. Based on an image created by Jordan Pierce \cite{wiki:Henon}.}
    \label{fig:henon_bifurcation}
\end{figure}

As a first example, let us consider the H\'enon map \cite{henon76} defined by
\begin{equation*}
    \begin{aligned}
        x_{j+1} &= 1-\lambda x_j^2 + y_j,\\
        y_{j+1} &= \nu x_j,
    \end{aligned}
\end{equation*}
where $ \lambda $ and $ \nu $ are parameters. Here, we assume that $ \nu = 0.3 $ is fixed and $\lambda \in \Lambda = [1.2, 1.4]$ is an uncertain parameter. The bifurcation diagram in Figure~\ref{fig:henon_bifurcation} illustrates the dynamical behavior for this parameter regime.

\begin{figure}[htb]
    \begin{minipage}{0.49\textwidth}
        \centering
        \subfiguretitle{(a) $\ell = 6$}
        \includegraphics[width = \textwidth]{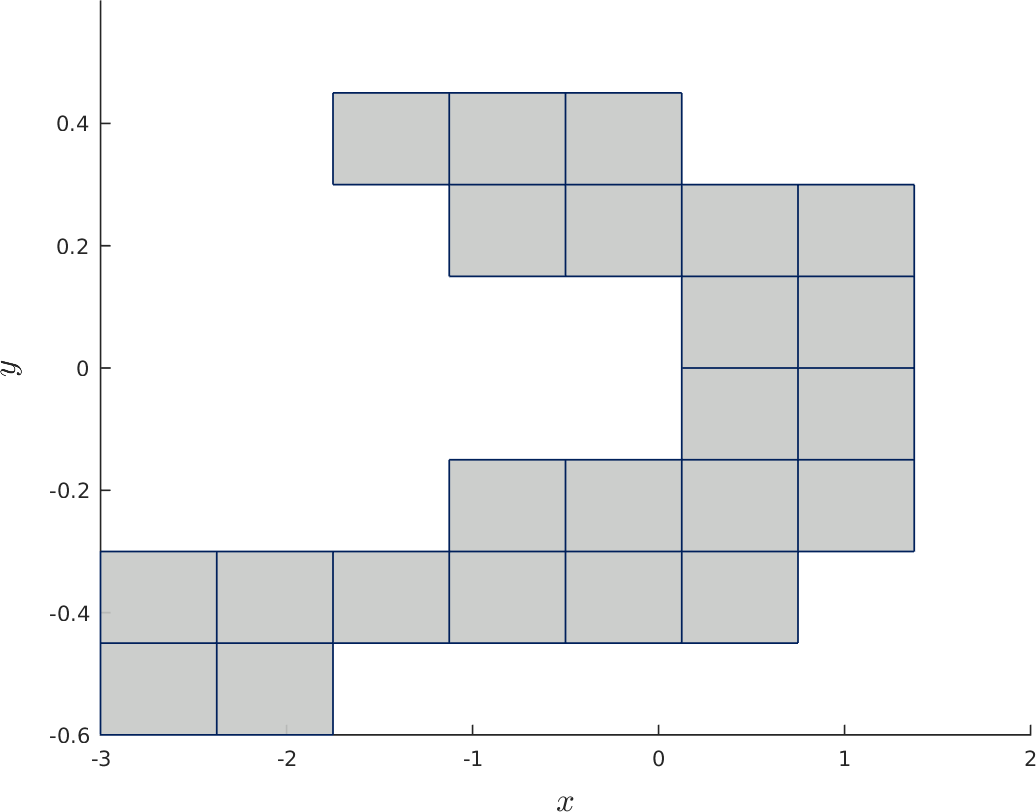}
    \end{minipage}
    \begin{minipage}{0.49\textwidth}
        \centering
        \subfiguretitle{(b) $\ell = 10$}
        \includegraphics[width = \textwidth]{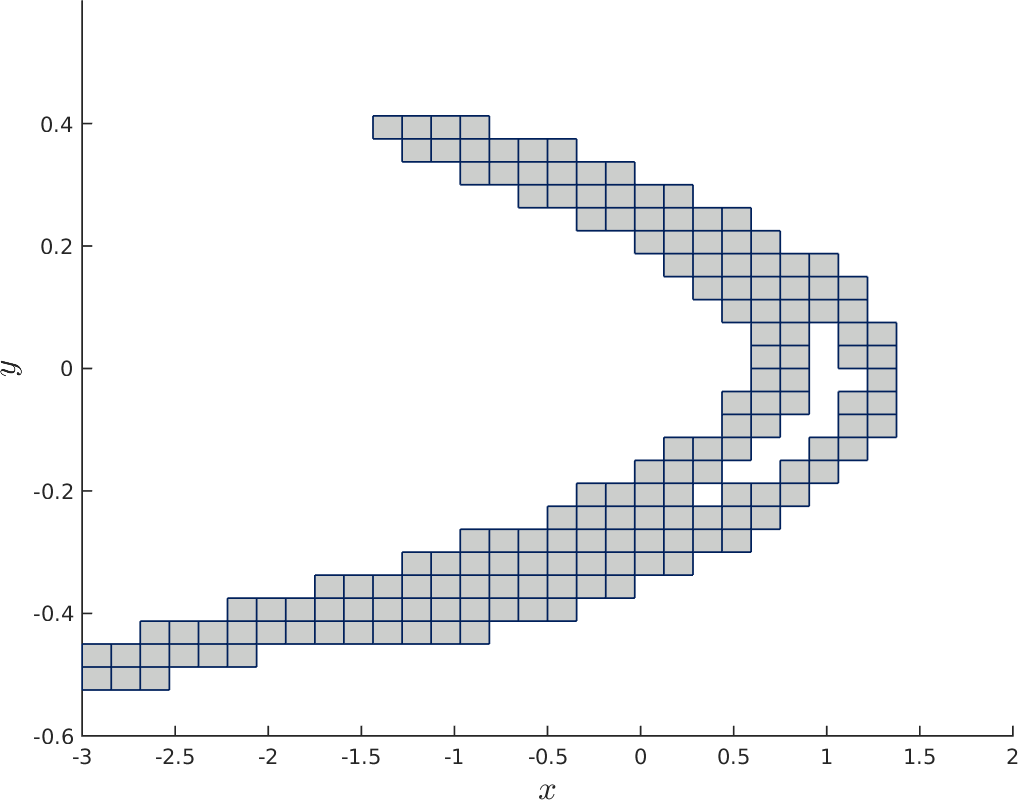}
    \end{minipage} \\[2ex]
    \begin{minipage}{0.49\textwidth}
        \centering
        \subfiguretitle{(c) $\ell = 14$}
        \includegraphics[width = \textwidth]{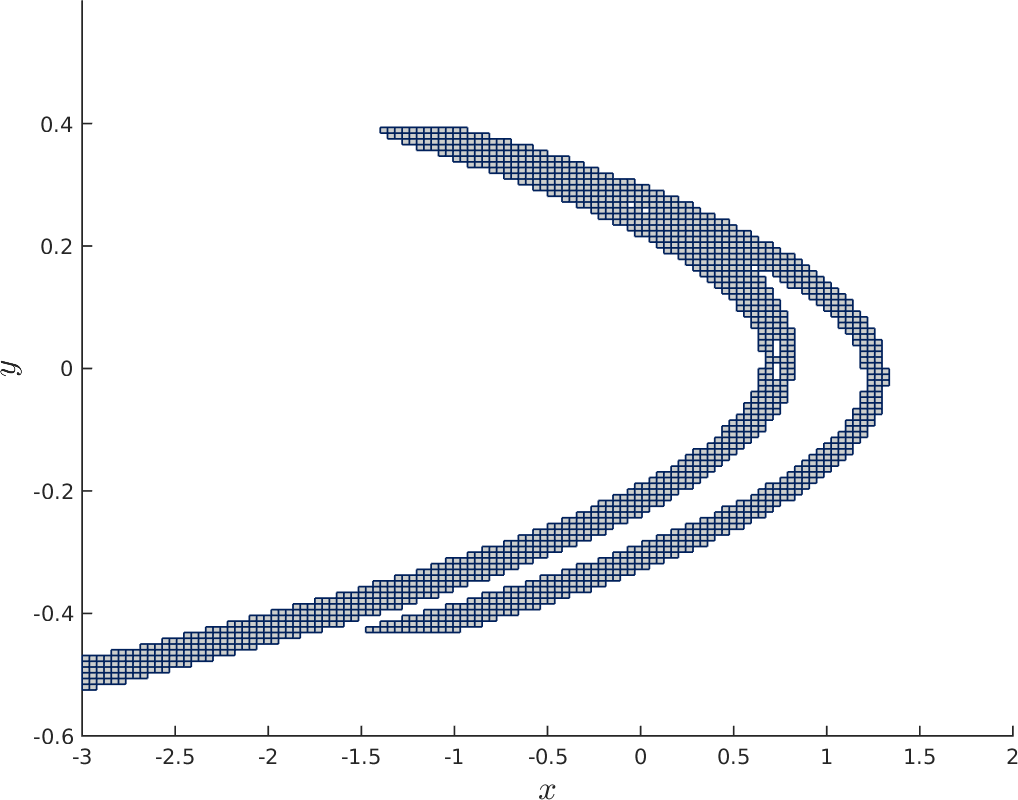}
    \end{minipage}
    \begin{minipage}{0.49\textwidth}
        \centering
        \subfiguretitle{(d) $\ell = 20$}
        \includegraphics[width = \textwidth]{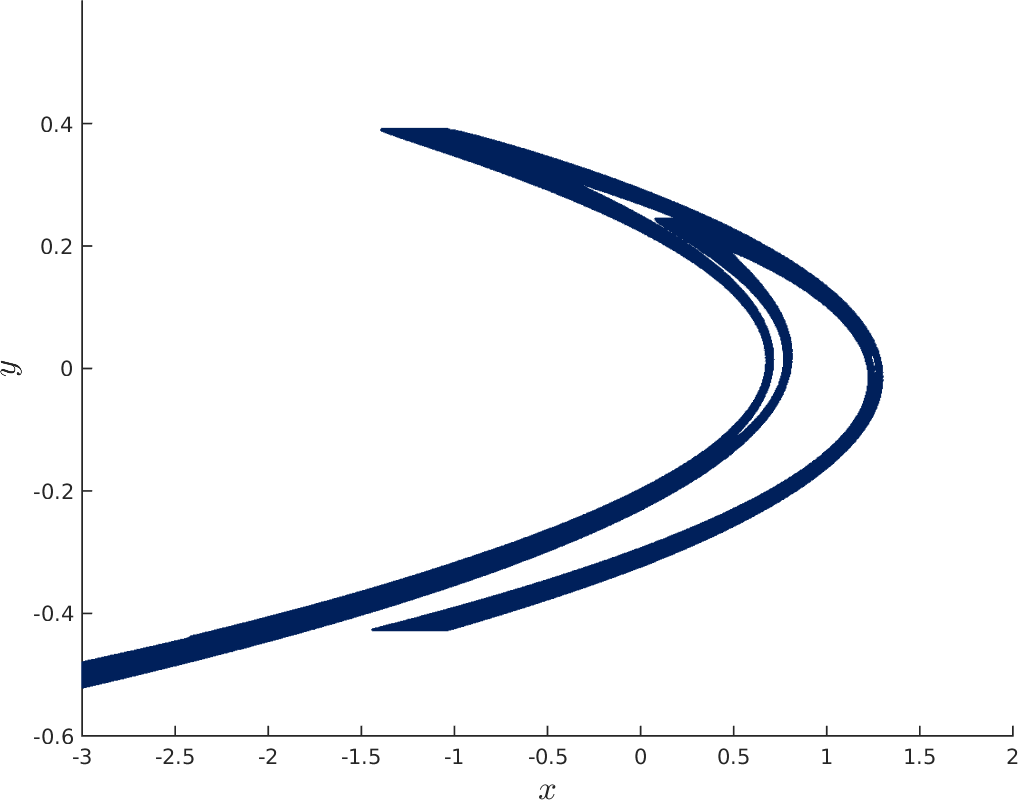}
    \end{minipage}
    \caption{Box coverings of the $(Q,\Lambda)$-attractor $A_{Q,\Lambda}$ of the H\'enon map obtained by the subdivision scheme after $ \ell $ subdivision steps.}
    \label{fig:henon_subdivision}
\end{figure}

We then approximate the $(Q,\Lambda)$-attractor $A_{Q,\Lambda}$ for $Q = [-3,2] \times [-0.6,0.6]$ using the subdivision algorithm described in Section~\ref{sec:comp_para_att}. In Figure~\ref{fig:henon_subdivision}, we show corresponding box coverings obtained by the algorithm after $ 6 $, $10 $, $14 $, and $ 20 $ subdivision steps.

\begin{figure}[htb]
    \begin{minipage}{0.49\textwidth}
        \centering
        \subfiguretitle{(a) $\lambda = 1.4$}
        \includegraphics[width = \textwidth]{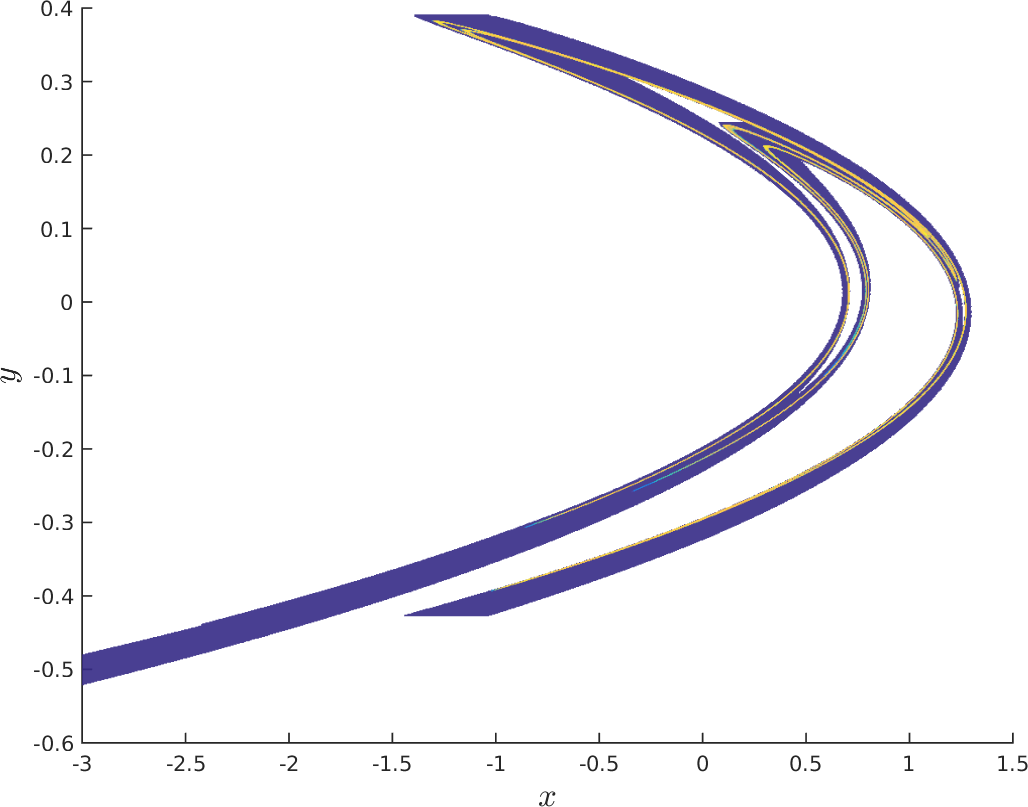}
    \end{minipage}
    \begin{minipage}{0.49\textwidth}
        \centering
        \subfiguretitle{(b) $\lambda = 1.24$}
        \includegraphics[width = \textwidth]{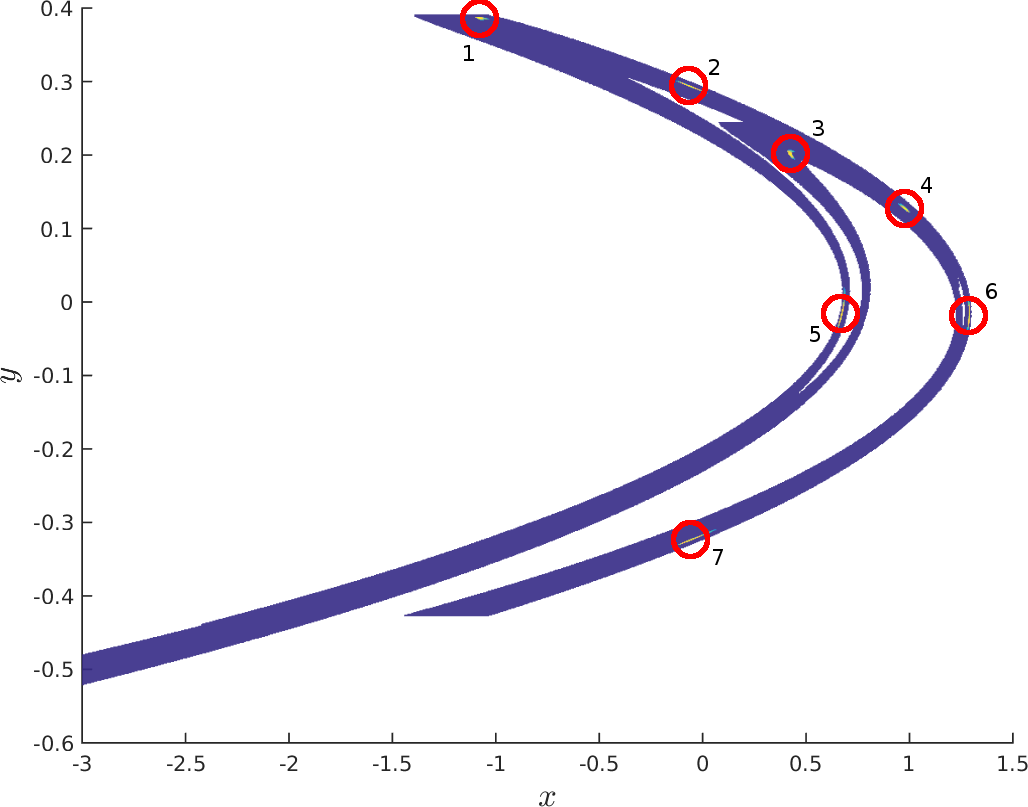}
    \end{minipage} \\[2ex]
    \begin{minipage}{0.49\textwidth}
        \centering
        \subfiguretitle{(c) $\lambda \sim \cN(1.24,0.0004)$}
        \includegraphics[width = \textwidth]{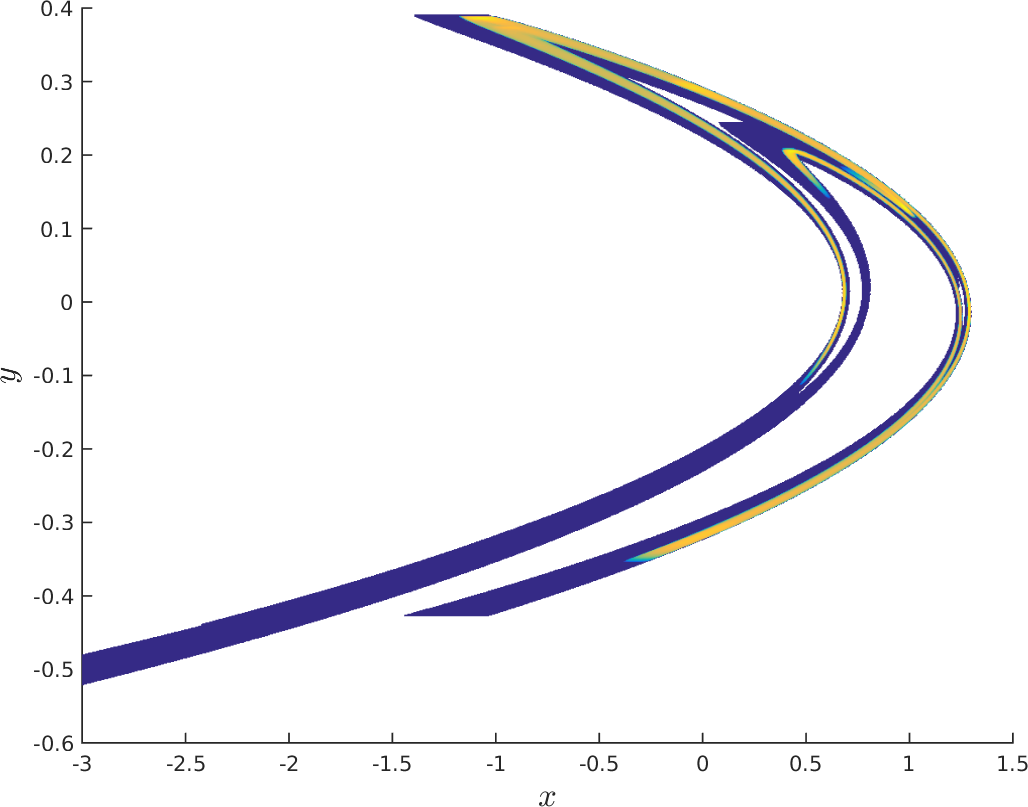}
    \end{minipage}
    \begin{minipage}{0.49\textwidth}
        \centering
        \subfiguretitle{(d) $\lambda \sim \cU(\Lambda)$}
        \includegraphics[width = \textwidth]{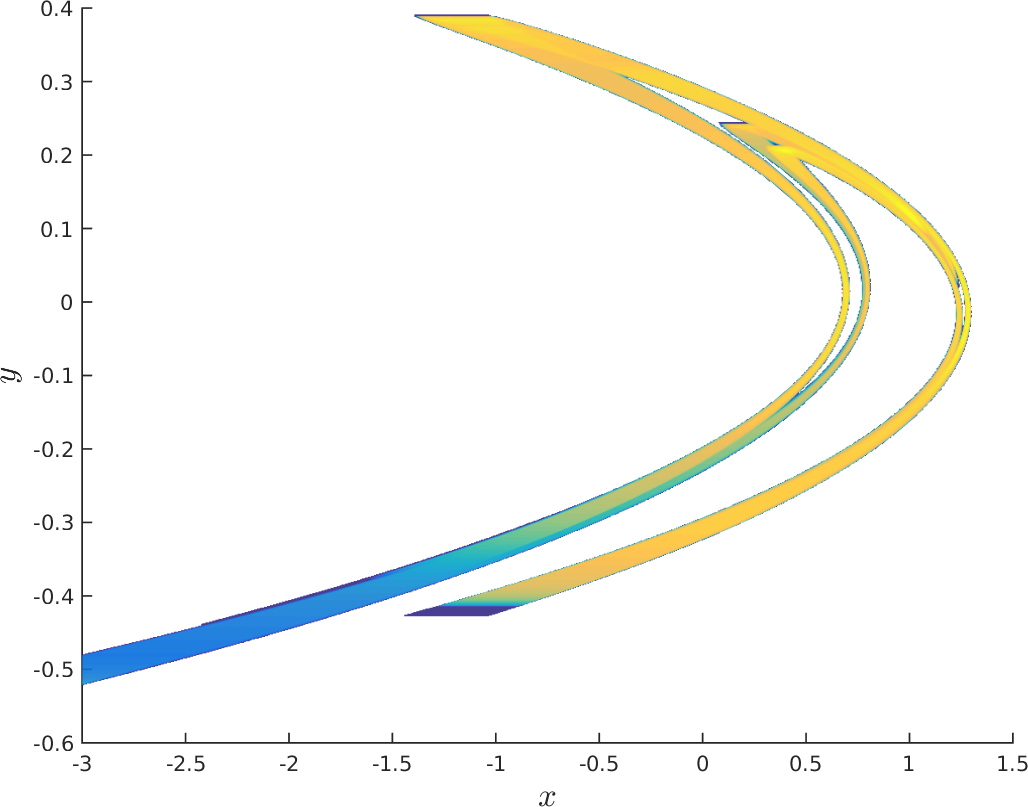}
    \end{minipage}
    \caption{Invariant measure on the $(Q,\Lambda)$-attractor $A_{Q,\Lambda}$ of the H\'enon map for different $\lambda$-distributions. The density ranges from blue (low density) $\rightarrow$ green $\rightarrow$ yellow (high density).}
    \label{fig:henon_measure}
\end{figure}

Given the box covering $\{C_1,\dots,C_d\}$ obtained by the subdivision scheme, we now
use Algorithm~\ref{alg2} for the approximation of invariant measures.
The resulting invariant measures for different $\lambda$-distributions
are shown in Figure~\ref{fig:henon_measure}.
Observe that the "seven-periodic behavior" in part (b) of the figure is still
visible in the result for the symmetrically truncated Gaussian in part (c).

\subsection{Van der Pol oscillator}
Let us now consider the van der Pol system, given by
\begin{equation}\label{eq:vdp}
    \begin{aligned}
        \dot x_1 &= x_2, \\
        \dot x_2 &= \lambda(1-x_1^2)x_2 - x_1,
    \end{aligned}
\end{equation}
where $ \lambda $ is the uncertain parameter. Here, we chose the interval $ \Lambda = [0.5,1.5] $. For each $\lambda \in \Lambda$ the system possesses a stable periodic solution as well as an unstable equilibrium in the origin. In our computation, we approximate the $(Q,\Lambda)$-attractor $A_{Q,\Lambda}$ for $Q = [-3,3]\times[-4,4]$. In this example, $f(x,\lambda)$ is given by the time-$T$-map $\phi^T(x,\lambda)$ with $T=4$, where $\phi$ is the flow of \eqref{eq:vdp}. We assume that after this time the full parameter uncertainty on $\Lambda$ is again relevant.

\begin{figure}[htb]
    \begin{minipage}{0.49\textwidth}
        \centering
        \subfiguretitle{(a) $\ell = 10$}
        \includegraphics[width = \textwidth]{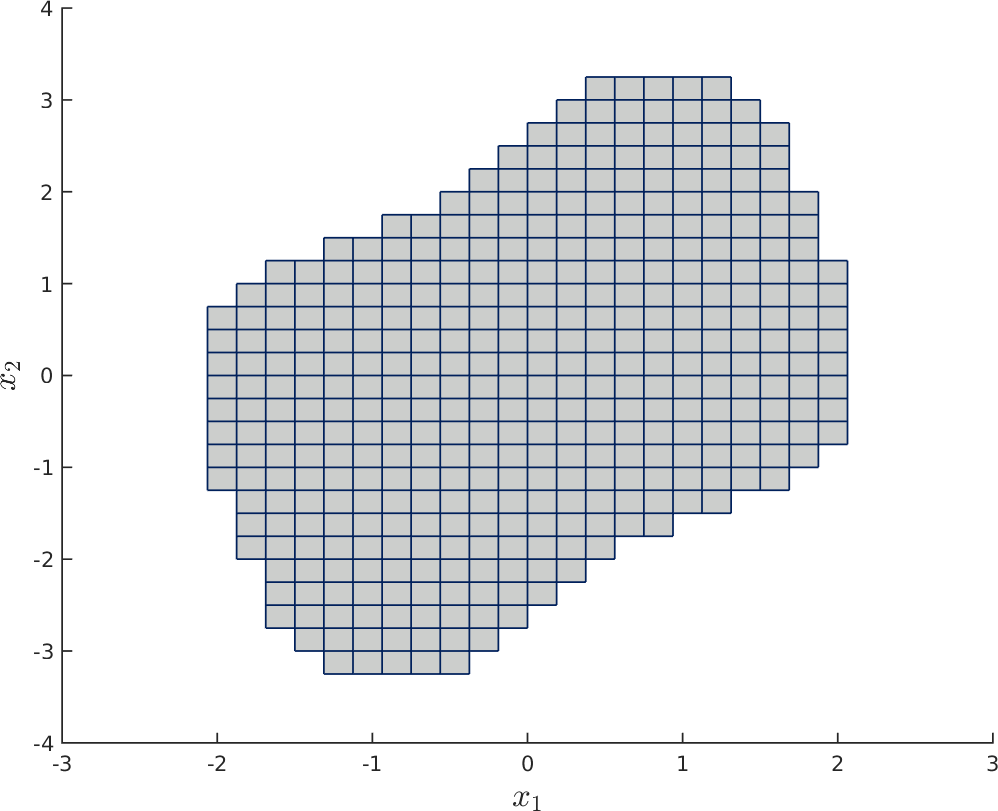}
    \end{minipage}
    \begin{minipage}{0.49\textwidth}
        \centering
        \subfiguretitle{(b) $\ell = 14$}
        \includegraphics[width = \textwidth]{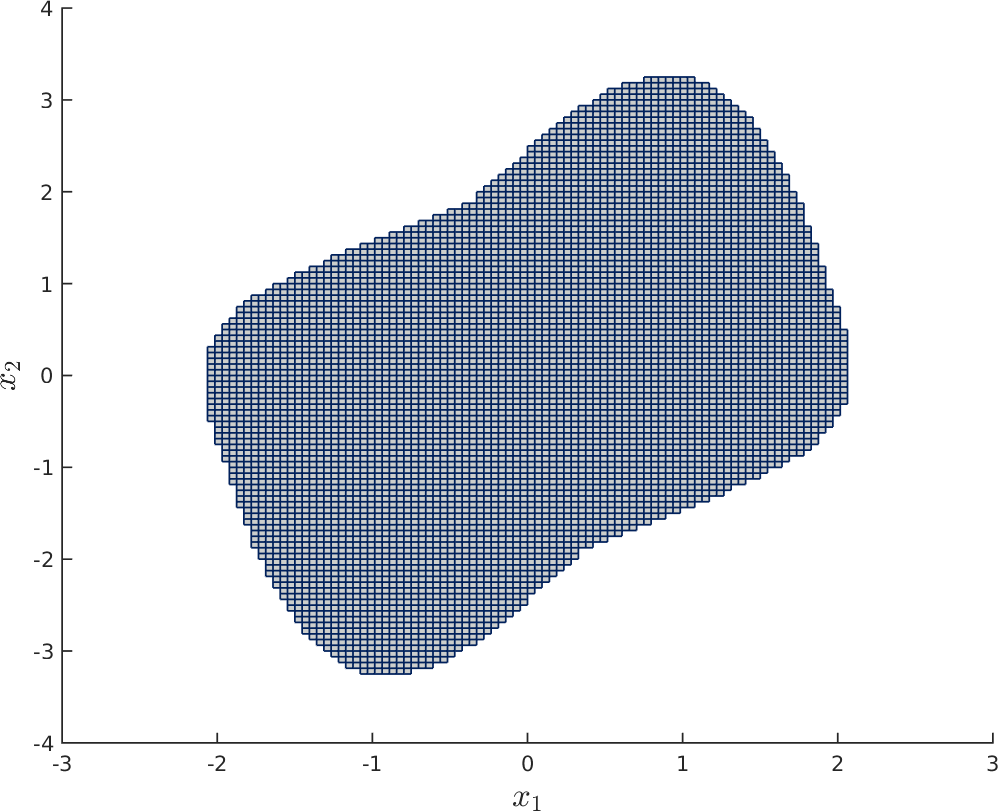}
    \end{minipage} \\[2ex]
    \begin{minipage}{0.49\textwidth}
        \centering
        \subfiguretitle{(c) $\ell = 18$}
        \includegraphics[width = \textwidth]{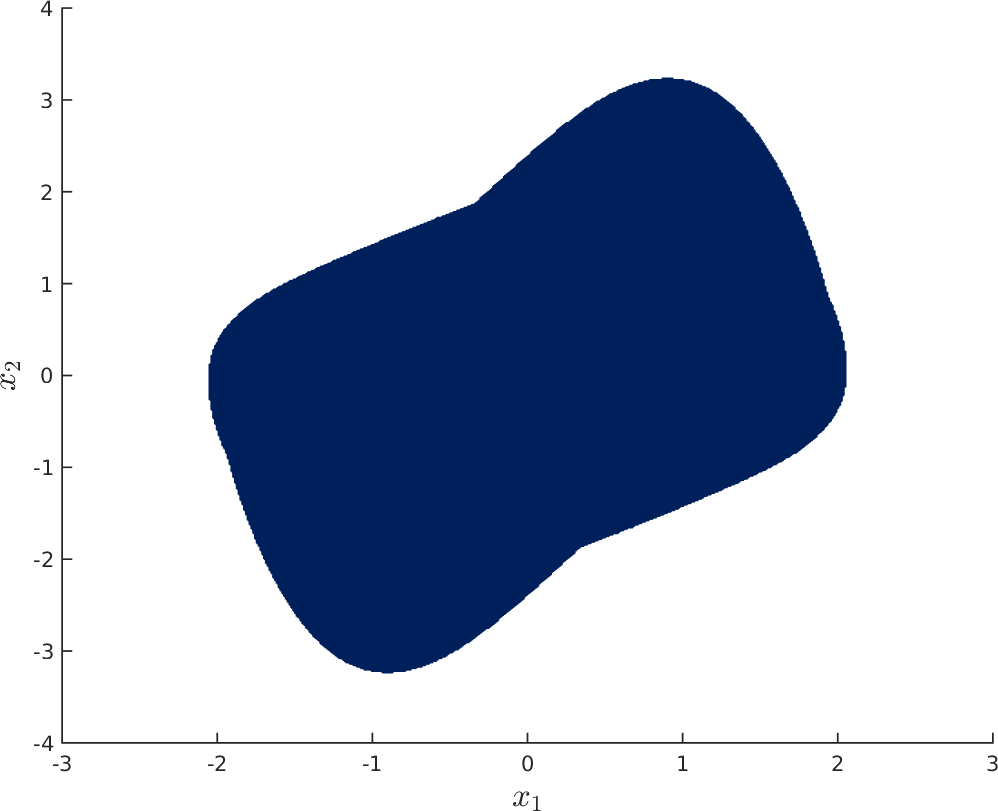}
    \end{minipage}
    \begin{minipage}{0.49\textwidth}
        \centering
        \subfiguretitle{(d) $\ell = 18$}
        \includegraphics[width = \textwidth]{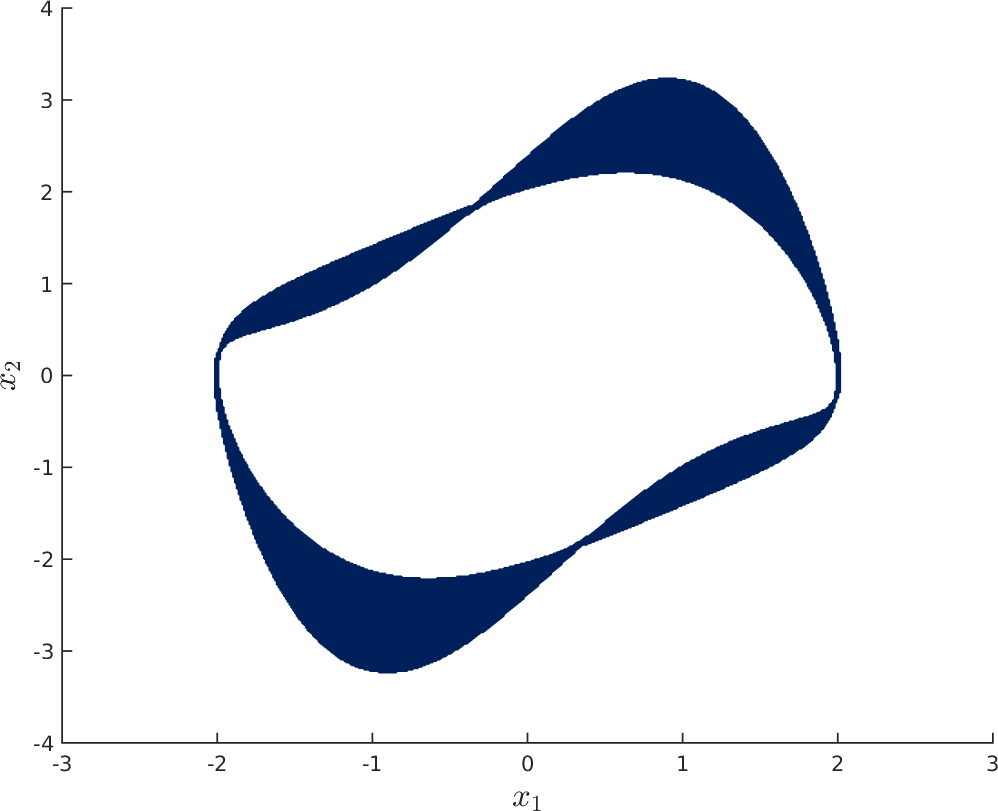}
    \end{minipage}
    \caption{Box coverings of the $(Q,\Lambda)$-attractor $A_{Q,\Lambda}$ of the van der Pol system obtained by the subdivision scheme after $ \ell $ subdivision steps. (a) -- (c) $Q = [-3,3]\times[-4,4]$. (d) $ Q $ without a small open neighborhood of the origin.}
    \label{fig:vdp_pu}
\end{figure}

Figure~\ref{fig:vdp_pu} (a)--(c) shows box coverings of the reconstruction of the two-dimensional unstable manifolds which accumulate on the stable periodic orbits at their boundary. Figure~\ref{fig:vdp_pu} (d) shows a box covering of the reconstructed periodic solutions itself. It has been obtained by removing a small open neighborhood $ U $ of the origin, i.e.~$ \widetilde{Q} = Q \backslash U $, resulting in the displayed box covering of $A_{\widetilde Q,\Lambda}$.

\begin{figure}[htb]
    \begin{minipage}{0.49\textwidth}
        \centering
        \subfiguretitle{(a) $\lambda = 1$}
        \includegraphics[width = \textwidth]{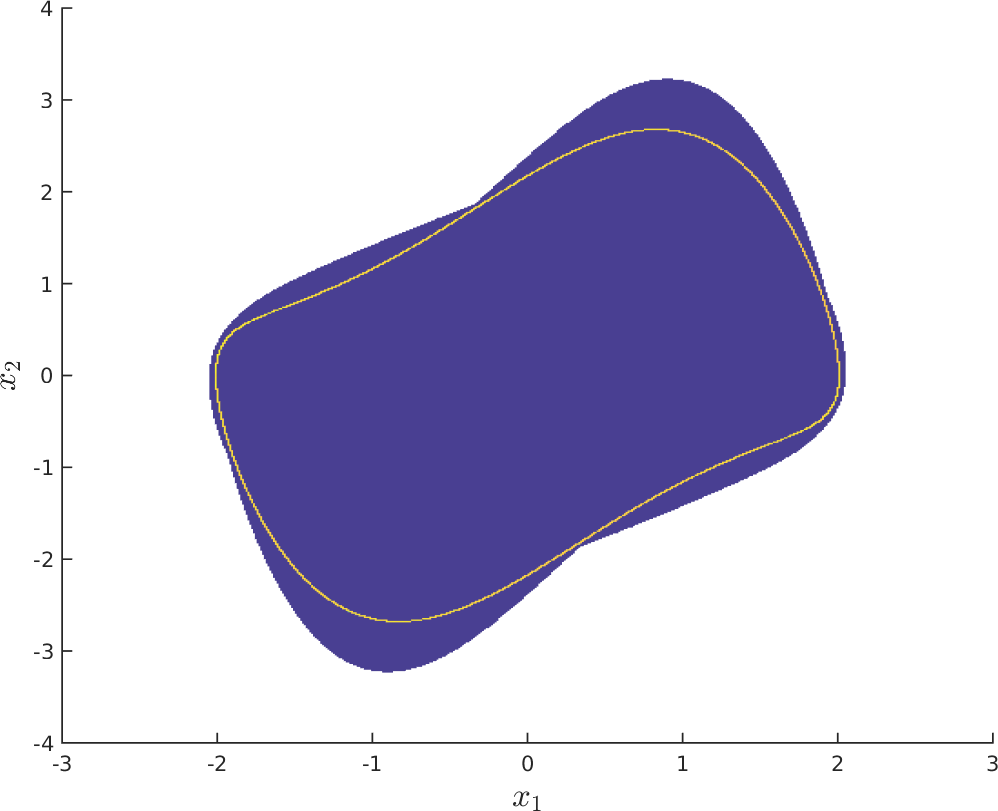}
    \end{minipage}
    \begin{minipage}{0.49\textwidth}
        \centering
        \subfiguretitle{(b) $\lambda \sim \cN(1,0.01)$}
        \includegraphics[width = \textwidth]{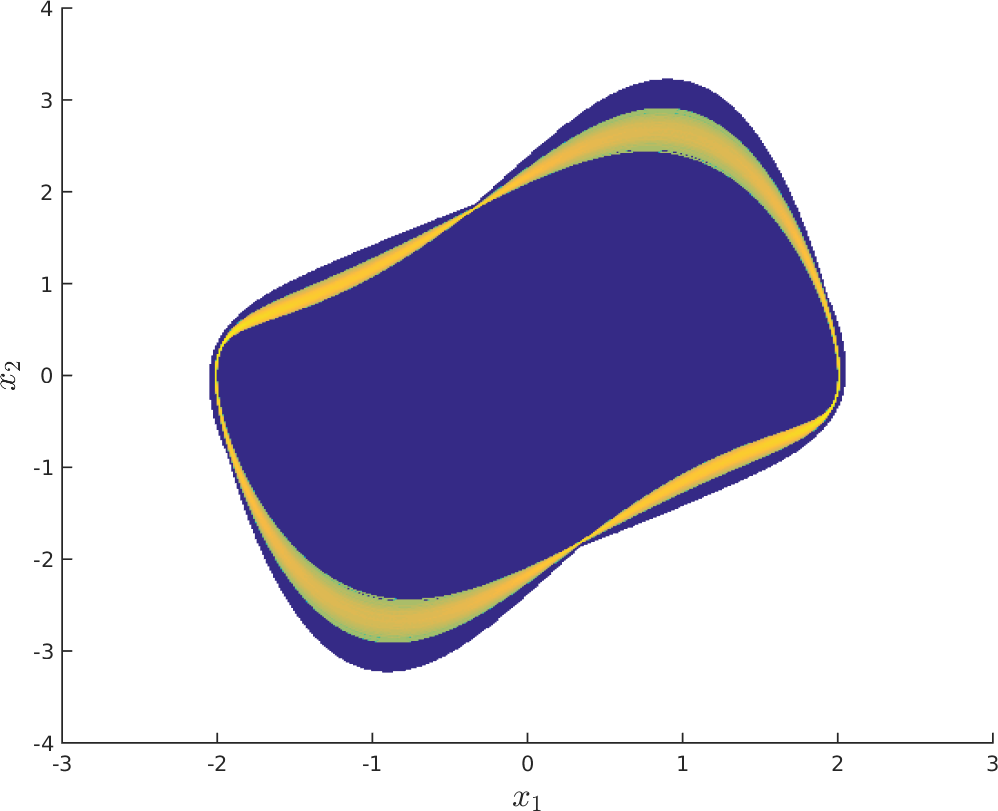}
    \end{minipage} \\[2ex]
    \begin{minipage}{0.49\textwidth}
        \centering
        \subfiguretitle{(c) $\lambda \sim \cN(1,0.04)$}
        \includegraphics[width = \textwidth]{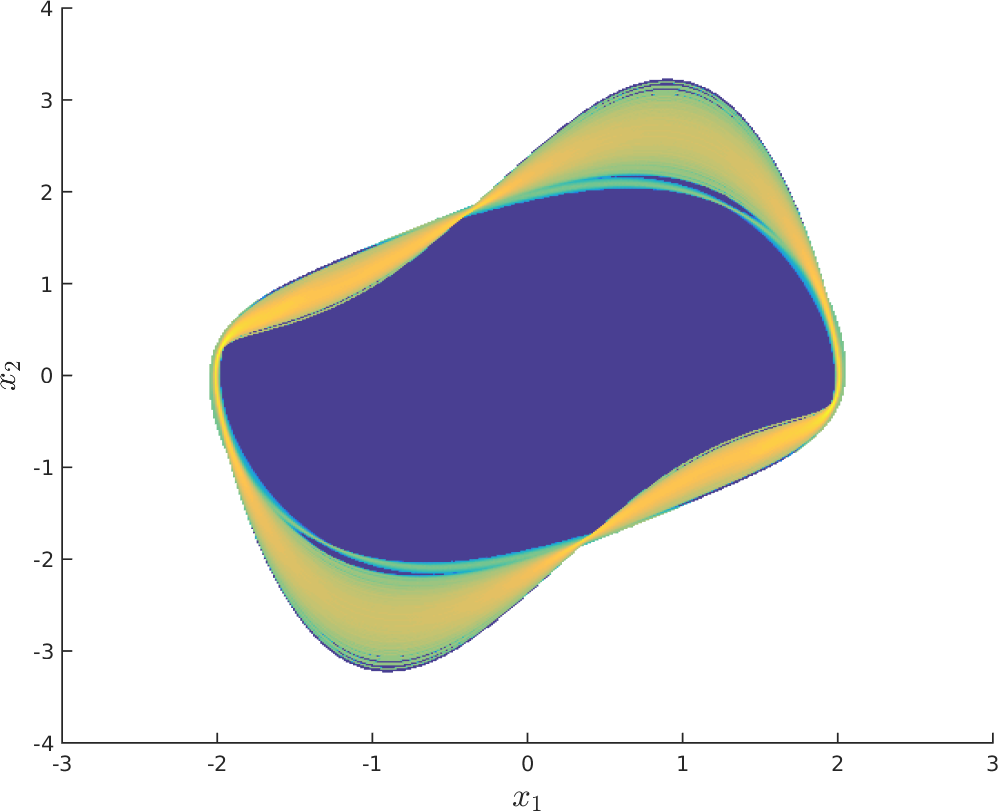}
    \end{minipage}
    \begin{minipage}{0.49\textwidth}
        \centering
        \subfiguretitle{(d) $\lambda \sim \cU(\Lambda)$.}
        \includegraphics[width = \textwidth]{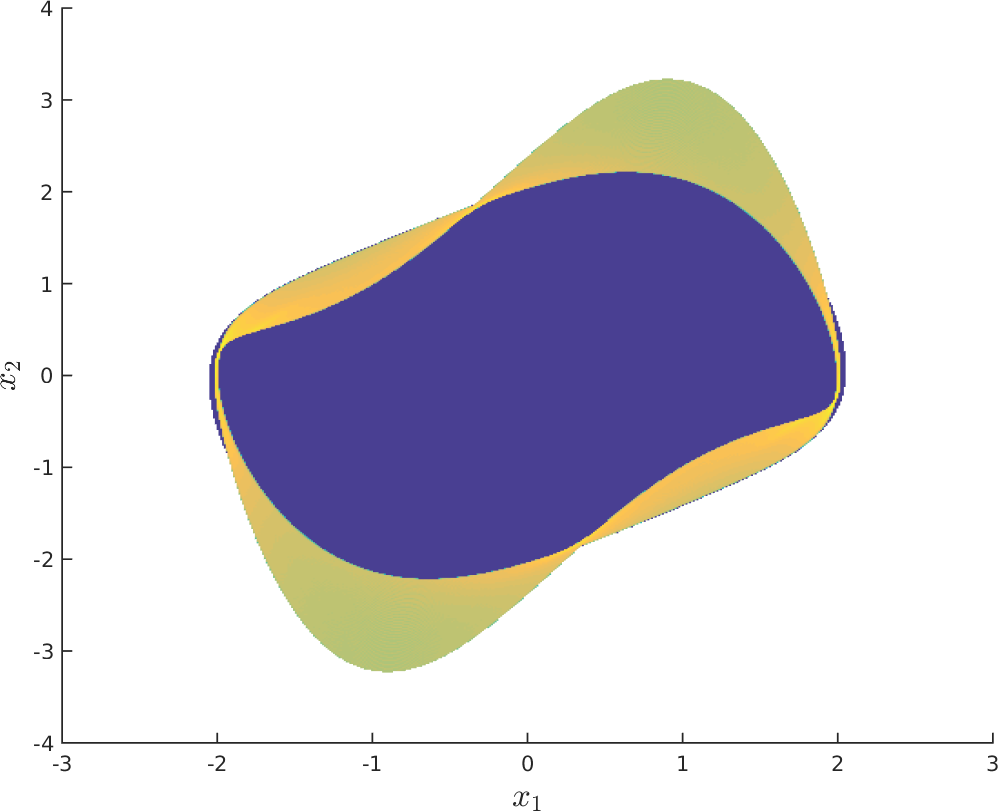}
    \end{minipage}
    \caption{Invariant measure on the $(Q,\Lambda)$-attractor $A_{Q,\Lambda}$ of the van der Pol system for different $\lambda$-distributions. The density ranges from blue (low density) $\rightarrow$ green $\rightarrow$ yellow (high density).}
    \label{fig:vdp_measure}
\end{figure}

Suppose that $\lambda \sim \cN(1,\sigma^2)$. In the case where $ \sigma = 0 $ (i.e.~$\lambda = 1$), we only have one stable periodic solution. In Figure \ref{fig:vdp_measure} (a), we show the corresponding invariant measure. Figure~\ref{fig:vdp_measure} (b)\&(c) shows the invariant measure for $\sigma = 0.1$ and $\sigma = 0.2$, Figure~\ref{fig:vdp_measure} (d) the corresponding measure assuming that $\lambda$ is uniformly distributed over $ \Lambda $.

\subsection{Arneodo system}

\begin{figure}[htb]
    \centering
    \includegraphics[width = 0.68\textwidth]{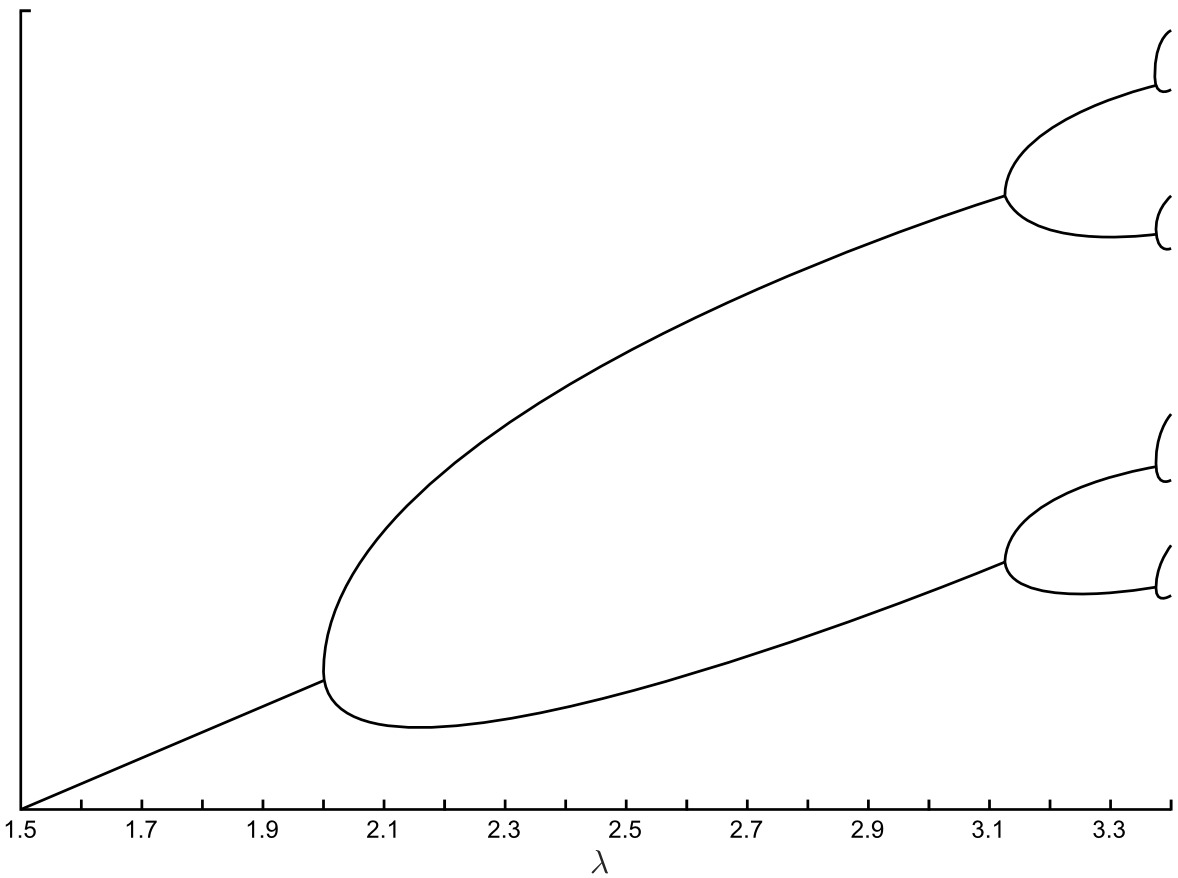}
    \caption{Schematic bifurcation diagram for the
    Arneodo system for $\lambda \in [1.5,3.4]$ (Hopf bifurcation and beginning of period
    doubling sequence).}
    \label{fig:arneodo_bifurcation}
\end{figure}

The final example is the Arneodo system~\cite{Arneodo1982}, which is given by
\begin{equation*}
    \frac{d^3 x}{dt^3} + \frac{d^2 x}{dt^2} + 2\frac{dx}{dt}-\lambda x+x^2 = 0.
\end{equation*}
We use the equivalent reformulation as a first-order system
\begin{equation}\label{eq:arneodo}
    \begin{aligned}
        \dot x_1 &= x_2,\\
        \dot x_2 &= x_3,\\
        \dot x_3 &= -x_3 -2x_2+\lambda x_1 -x_1^2.
    \end{aligned}
\end{equation}
This system possesses the equilibria $\widetilde X_1 = (0,0,0)$ and $\widetilde X_2(\lambda) = (\lambda,0,0)$. The latter is asymptotically stable for $\lambda < 2$. At $\lambda = 2$, the equilibrium
$\widetilde X_2(\lambda)$ undergoes a supercritical Hopf bifurcation (cf.~\cite{KrOs1999}). For $\lambda > 2$,
points on the two-dimensional unstable manifold of $\widetilde X_2(\lambda)$ converge to the corresponding
limit cycle on the branch of periodic solutions, where the amplitude of the limit cycle
grows with increasing values of $\lambda$. In Figure \ref{fig:arneodo_bifurcation},
we show a bifurcation diagram for the periodic solution for $\lambda \in [1.5,3.4]$.
For $\lambda \approx 3.1$, the limit cycle loses its stability in a period-doubling bifurcation. We 
choose $\Lambda = [2.8,3.4]$ in order to quantify the uncertainty in this region where the system undergoes several bifurcations. Analogous to the second example, $f(x,\lambda)$ is given by the time-$T$-map $\phi^T(x,\lambda)$ with $T=2$, where $\phi$ is the flow of \eqref{eq:arneodo}. We assume that after this time the full parameter uncertainty on $\Lambda$ comes again into play.

\begin{figure}[htb]
    \begin{minipage}{0.49\textwidth}
        \centering
        \subfiguretitle{(a) $\ell = 12$}
        \includegraphics[width = \textwidth]{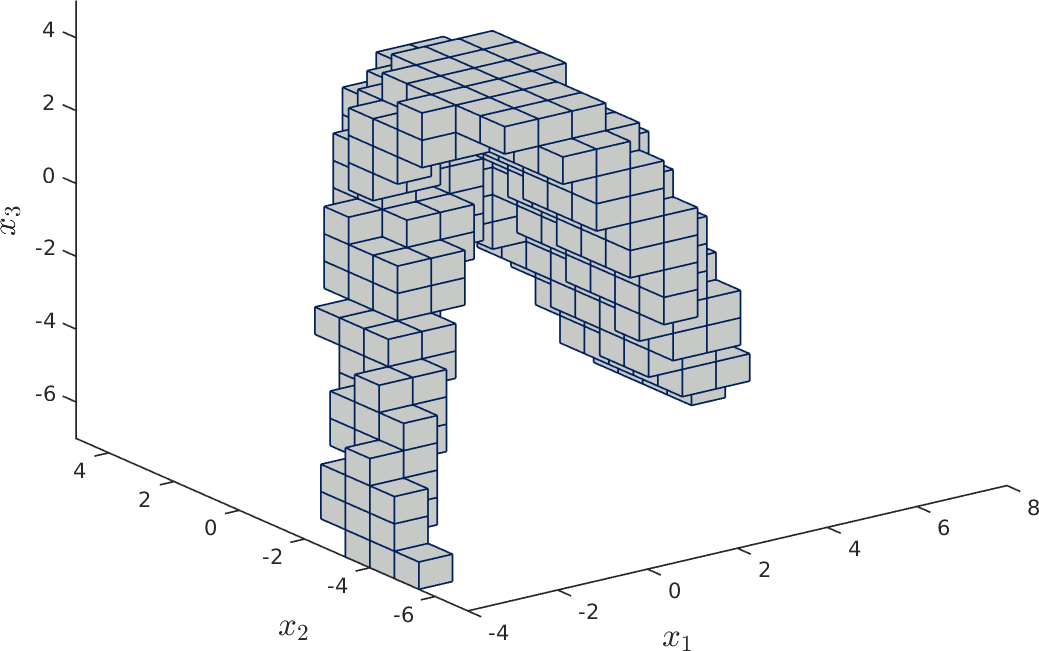}
    \end{minipage}
    \begin{minipage}{0.49\textwidth}
        \centering
        \subfiguretitle{(b) $\ell = 18$}
        \includegraphics[width = \textwidth]{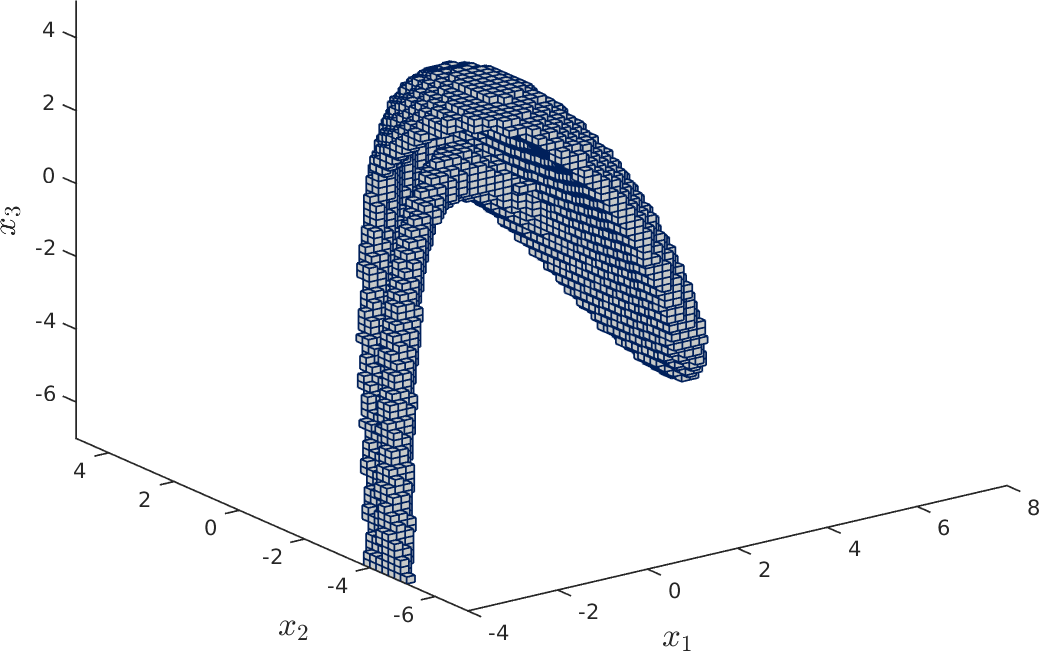}
    \end{minipage} \\[2ex]
    \begin{minipage}{0.49\textwidth}
        \centering
        \subfiguretitle{(c) $\ell = 24$}
        \includegraphics[width = \textwidth]{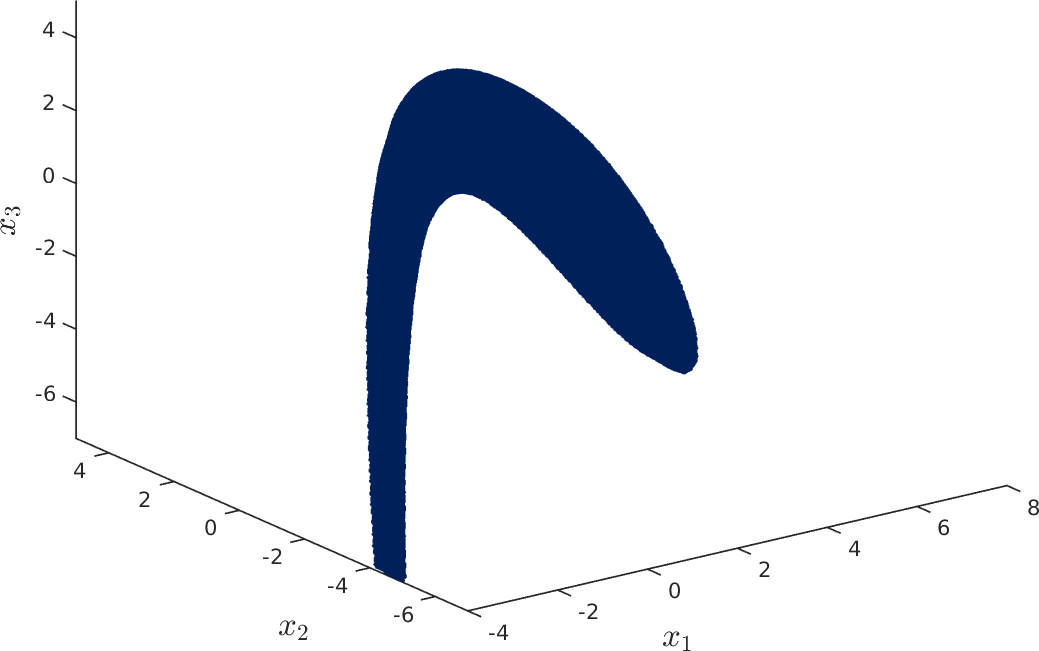}
    \end{minipage}
    \begin{minipage}{0.49\textwidth}
        \centering
        \subfiguretitle{(d) $\ell = 24$}
        \includegraphics[width = \textwidth]{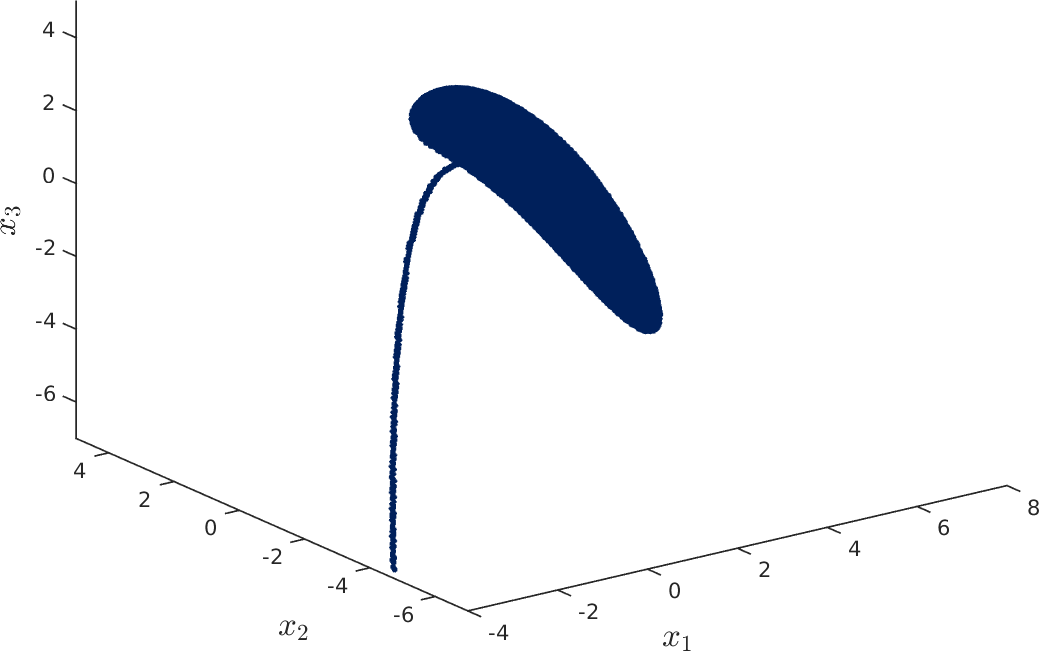}
    \end{minipage}
    \caption{(a)--(c) Box coverings of the $(Q,\Lambda)$-attractor $A_{Q,\Lambda}$ of the Arneodo system obtained by the subdivision scheme after $ \ell $ subdivision steps, $\Lambda = [2.8,3.4]$.
    (d)  $A_{Q,\Lambda}$ for fixed $\bar{\lambda} = 3.1$.}
    \label{fig:arneodo_pu}
\end{figure}

\begin{figure}[htb]
    \begin{minipage}{0.49\textwidth}
        \centering
        \includegraphics[width=\textwidth]{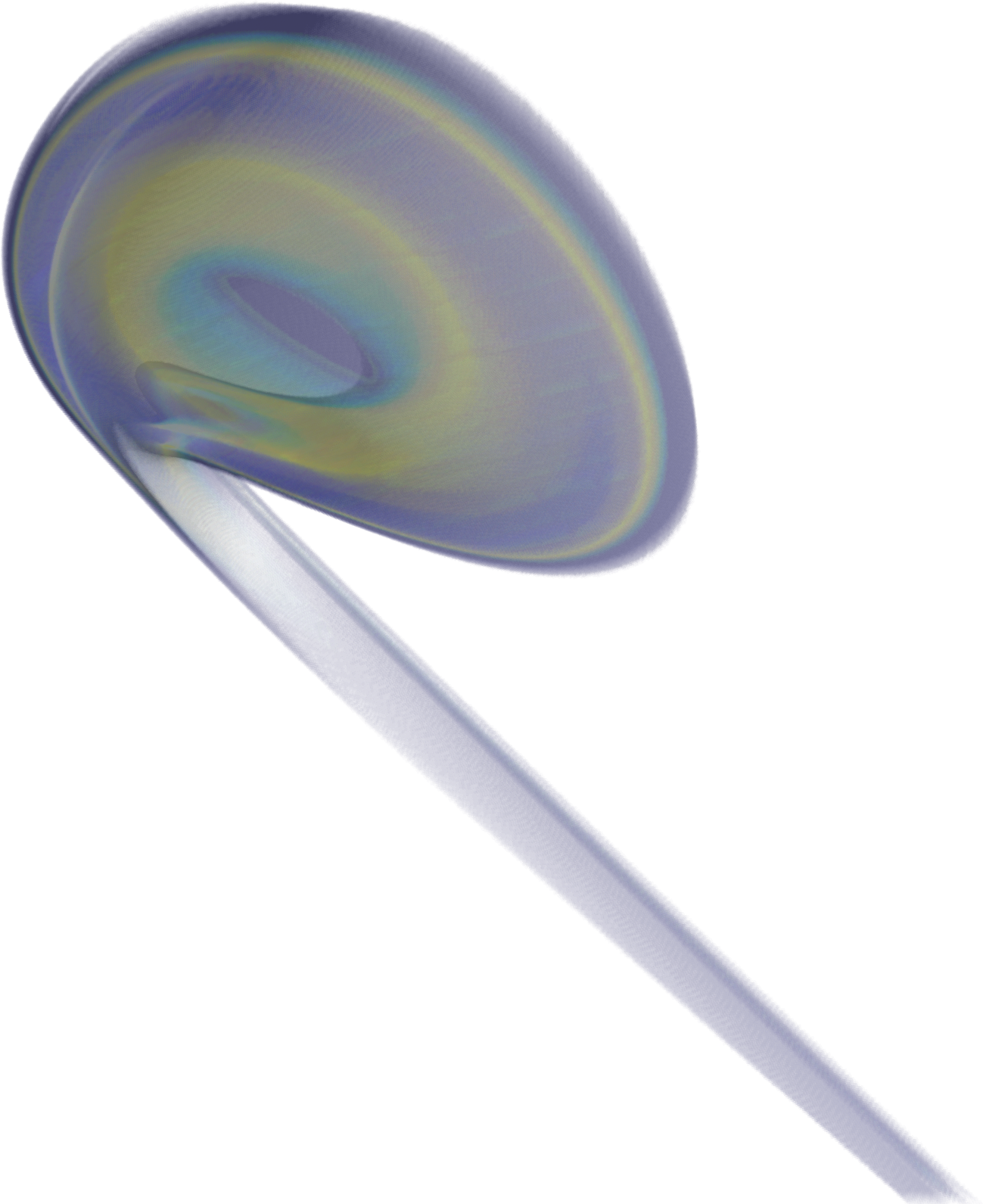}
    \end{minipage}
    \hfill
    \begin{minipage}{0.49\textwidth}
        \centering
        \includegraphics[width=\textwidth]{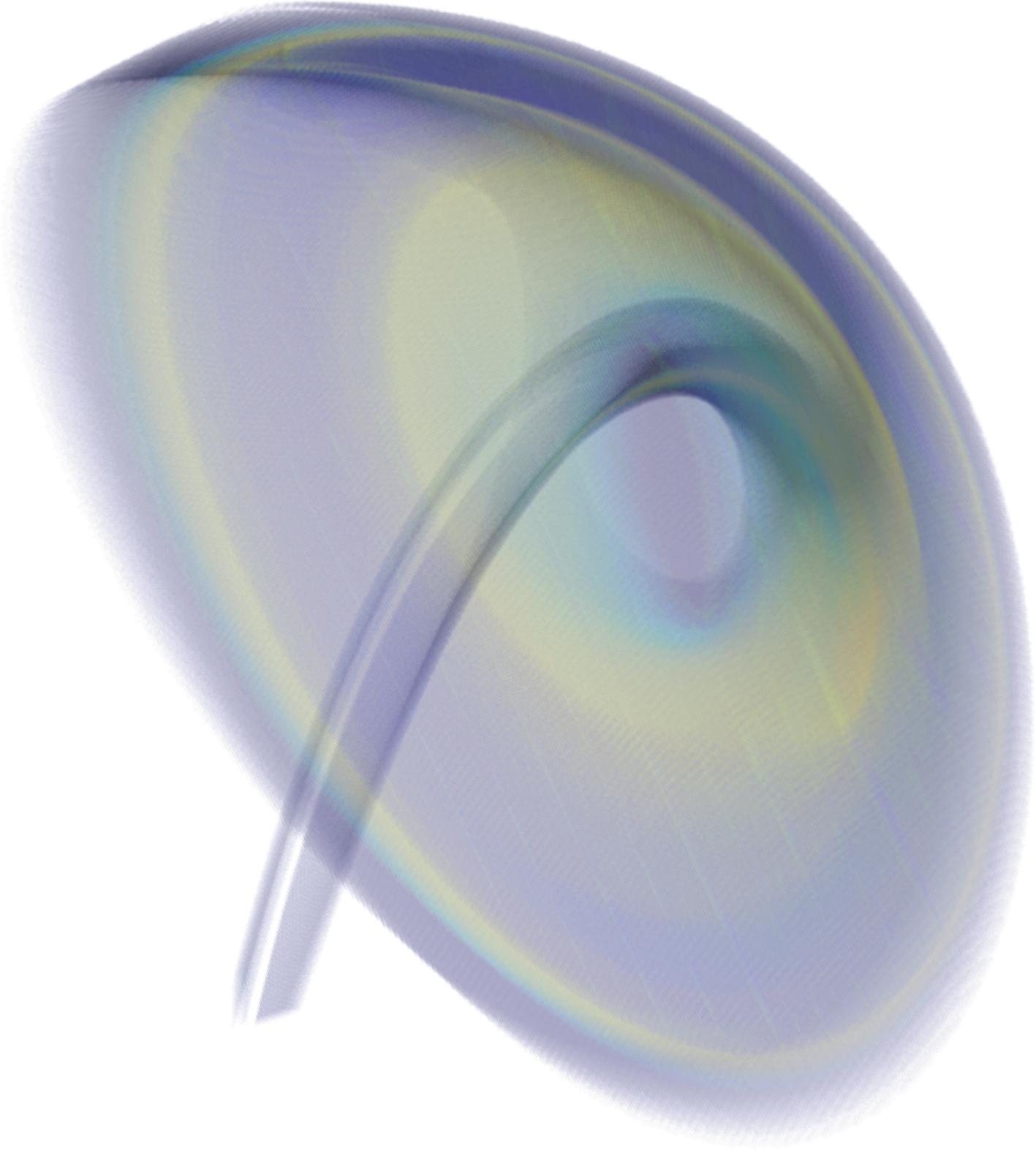}
    \end{minipage}
    \caption{Two projections of the invariant measure on the $(Q,\Lambda)$-attractor $A_{Q,\Lambda}$ of the Arneodo system for $\lambda \sim \cN(3.1,0.01)$.The density ranges from blue (low density) $\rightarrow$ green $\rightarrow$ yellow (high density).}
    \label{fig:arneodo_measure}
\end{figure}

Figure~\ref{fig:arneodo_pu} (a)--(c) shows successively finer box coverings of the $(Q,\Lambda)$-attractor $A_{Q,\Lambda}$ for $Q = [-4,8]\times[-7,5]\times[-7,5]$. In this way, we compute a reconstruction of two-dimensional unstable manifolds of $\widetilde X_2(\lambda)$ which either accumulate on limit cycles, the period-doubled limit cycles, or even higher periodic limit cycles, depending on the value of $\lambda$. We also obtain a covering of the one-dimensional unstable manifold of $\widetilde X_1$. For comparison purposes, Figure \ref{fig:arneodo_pu} (d) depicts the $(Q,\Lambda)$-attractor $A_{Q,\Lambda}$ for a fixed value of $\bar{\lambda} = 3.1$, i.e.~without an underlying parameter uncertainty. In Figure~\ref{fig:arneodo_measure}, we finally show two projections of the invariant measure on $A_{Q,\Lambda}$, where $\lambda$ is a Gaussian.

\section{Conclusion}
\label{sec:Conclusion}

In this paper, we introduce the notion of $(Q,\Lambda)$-attractors, which can be regarded as a generalization of relative global attractors for dynamical systems with uncertain parameters. These objects capture all the dynamics
which may potentially be induced by a given parameter distribution on $\Lambda$. We then use a classical
transfer-operator approach to compute invariant measures on $(Q,\Lambda)$-attractors
which are related to different $\lambda$-distributions. The technical framework from \cite{DJ99}
allows us to obtain a related convergence result in the context of small random perturbations.
The numerical examples illustrate the fact that the techniques are very well applicable to
low dimensional dynamical systems with parameter uncertainty.

So far, we considered only dynamical systems with one uncertain parameter, i.e.~$\Lambda \subset \R $. Future work
will include analyzing systems with multiple uncertain parameters and also the scalability of the proposed algorithms.
Due to the curse of dimensionality, analyzing high-dimensional systems or systems with a large number of uncertain parameters is in general challenging. Finally a set oriented numerical method adapted to the analysis of
uncertainty quantification with respect to initial conditions will also be developed.

\bibliographystyle{alpha}
\bibliography{UQ}

\end{document}